\input amstex
\documentstyle{amsppt}
\magnification=\magstep1
\input epsf.tex

\newcount\sectionno\sectionno=0
\newcount\equationno\equationno=0
\newcount\theoremno\theoremno=0

\font\bigfont=cmbx10
\def\Heading#1{\bigskip\bigskip\goodbreak\centerline{\bigfont #1}\nobreak\bigskip\nobreak}
\def\Section#1{\global\advance\sectionno by 1\relax%
   \equationno=0\theoremno=0\Heading{\the\sectionno. #1}}

\newwrite\aux
\immediate\openout\aux=\jobname.aux

\def\formeq{\the\sectionno.\the\equationno}
\def\formth{\the\sectionno.\the\theoremno}

\def\Elabel#1{\global\advance\equationno by 1 \formeq%
\immediate\write\aux{\def\string#1{\formeq}}%
\global\edef#1{\formeq}}

\def\Tlabel#1{\global\advance\theoremno by 1 \formth%
\immediate\write\aux{\def\string#1{\formth}}%
\global\edef#1{\formth}}

\def\gk{\kappa}
\def\gl{\lambda}
\def\x{\xi}
\def\gep{\varepsilon}
\def\ga{\alpha}
\def\gb{\beta}
\def\gga{\gamma}
\def\gG{\Gamma}
\def\gD{\Delta}
\def\gs{\sigma}
\def\gd{\delta}
\def\gth{\theta}
\def\dilpar{q}
\def\Bern{R}
\def\Energy{E}

\def\len{\text{len}}
\def\BbR{\Bbb R}
\def\BbC{\Bbb C}

\def\splice{\sqcup}
\def\T#1{{\vec #1}}
\def\norm#1{{\Vert #1 \Vert}}
\def\set#1{\{#1\}}
\def\abs#1{\left\vert {#1} \right\vert}
\def\pr#1{\left( {#1} \right)}
\def\iMATH{ i}
\def\IM{\mathop{\roman Im}}
\def\RE{\mathop{\roman Re}}
\def\diam{\mathop{\roman diam}}
\def\wt{\widetilde}
\def\dd{\,{\roman d}}
\def\BirkBoor{1}
\def\BirkBurch{2}
\def\Brunnett{3}

\def\Goss{5}
\def\GolombJerome{4}
\def\Jeromeone{6}
\def\LeeForsythe{7}
\def\Linner{8}
\def\Royden{9}

\def\Uturn{5.5}
\def\leftright{5.6}
\def\Gderiv{5.11}
\def\mostofit{5.12}
\def\MajorTom{6.2}
\def\bridge{6.8}

\hsize=6.5truein
\vsize=9truein
\overfullrule=0pt
\topmatter
\title
Elastic Splines I: Existence
\endtitle
\affil Dept. of Mathematics, Faculty of Science \\ Kuwait University \endaffil
\abstract Given interpolation points $P_1,P_2,\; \ldots,P_m$ in the plane, it is known
that there does not exist an interpolating curve with minimal bending energy, unless the given points lie sequentially
along a line. We say that an interpolating curve is {\it admissible} if each piece, connecting two consecutive points 
$P_i$ and $P_{i+1}$, is an
s-curve, where an {\it s-curve} is a 
planar curve which first turns monotonically at most $180^\circ$ in one direction and then turns monotonically at most $180^\circ$
in the opposite direction. 
Our main result is that among all admissible interpolating curves there exists a curve with minimal bending energy.
We also prove, in a very constructive manner, the existence of an s-curve, with minimal bending energy, that connects two given
unit tangent vectors.
\endabstract

\author Albert Borb\'ely \& Michael J. Johnson  \endauthor
    \address Department of Mathematics, Faculty of
    Science, Kuwait University, P.O. Box 5969, Safat 13060, Kuwait
    \endaddress
\date January 21, 2014 \enddate
\email borbely\@sci.kuniv.edu.kw, yohnson1963\@hotmail.com \endemail
\thanks Cite this article as:\newline
    Borb\'ely, A. \& Johnson, M.J. Constr Approx (2014) 40: 189--218. \hskip1truecm doi: 10.1007/s00365-014-9244-4\newline
The final publication is available at http://link.springer.com/article/10.1007/s00365-014-9244-4
\endthanks
 
\endtopmatter

\document

\Section{Introduction}
Given a sequence of points $P_1,P_2,\ldots,P_m$ in $\BbR^2$ with $P_i \neq P_{i+1}$, a curve $F:[a,b]\to \BbR^2$
is called an {\bf interpolating curve} if there exist times $a=t_1<t_2<\cdots < t_m=b$ such that $F(t_i)=P_i$.
In the special case when the interpolation points can be written as $P_i=(x_i,y_i)$, with $x_1<x_2<\cdots<x_m$, 
an interpolating curve can be constructed as the graph of a smooth function $g:[x_1,x_m]\to\BbR$, provided
$g$ satisfies the interpolation conditions $g(x_i)=y_i$. It is well known that if $g$ is the natural
cubic spline, then $g$ minimizes the functional
$\int_{a}^{b} (g''(x))^2\dd x$ among all smooth functions which satisfy the interpolation conditions. This functional
is often viewed as a simple approximation of the curve's bending energy
$\int_0^L \gk(s)^2\dd s$, where $s$ denotes arclength and $\gk$ denotes signed curvature,
and it is natural to ask what would happen if one tried to minimize the bending
energy among all smooth interpolating curves.  Unfortunately, such optimal curves do not exist except
in the trivial case  when the interpolation points lie sequentially along a line. Apparently,
this was first observed by Birkhoff and de Boor \cite{\BirkBoor}, along with Birkhoff, Burchard and 
Thomas \cite{\BirkBurch}. This lack of existence can be understood as a consequence of the effect that scaling has on
bending energy: the bending energy of a curve scaled by a factor $\dilpar$ equals $\frac1{\dilpar}$ times the
original bending energy. As a result, it is possible to construct smooth interpolating curves with
arbitrarily small bending energy. For example, let $p\gg 1$ and consider the circles $c_1,c_2,\ldots,c_m$ which meet 
tangentially at the point $(p,0)$; specifically, let
$c_i$ be the circle which begins and ends at $(p,0)$, has center
on the $x$-axis, and passes through the point $P_i$. Then the subcurve of $c=c_1\cup c_2 \cup \cdots \cup c_m$, starting
at $P_1$ and ending at $P_m$, is an interpolating curve
whose bending energy tends to $0$ as $p\to\infty$
(see section 2.6 of \cite{\Linner} for other constructions).

Subsequent attention
was directed towards interpolating curves whose bending energy is locally minimal (i.e. minimal among
all nearby interpolating curves).  It was reported in \cite{\BirkBurch}, and mentioned in
\cite{\BirkBoor}, that
if an interpolating curve $F$ has a locally minimal bending energy, then each segment of $F$, connecting two consecutive
interpolation points, will be a segment of `rectangular elastica', meaning a planar curve
whose signed curvature $\gk$ satisfies the differential equation $2\frac{d^2\gk}{ds^2} +\gk^3=0$.
(Rectangular elastica was first described by James Bernoulli (1694) and is
one of the nine types of elastica identified by Euler (1750), see \cite{\Goss}.)
Using a variational calculus and physical reasoning, Lee and Forsythe \cite{\LeeForsythe} (see also 
\cite{\Brunnett}) have confirmed that each
segment of $F$ is indeed a segment of rectangular elastica, and have moreover shown that the signed
curvature of $F$ is continuous throughout the curve and vanishes at the endpoints. 
Unfortunately, interpolating curves with
locally minimal bending energy do not necessarily exist, and this constitutes a significant deficiency in the theoretical
foundation of this interpolation method. 

Rather than seeking an interpolating curve with a locally minimal bending energy, an alternate approach is to define
a restricted class of `admissible' interpolating curves and then seek a curve with minimal bending
energy in the restricted class. Birkhoff proposed a restriction on {\it length} and conjectured that among all
smooth interpolating curves of length at most $L_0$, $L_0$ being a prescribed upper bound, 
there exists a curve with minimal
bending energy.  This conjecture was eventually proved by Jerome \cite{\Jeromeone} (see \cite{\GolombJerome} for a more
comprehensive treatment and also \cite{\Linner} where `pinning' and `clamping' at interpolation nodes are treated).

Rather than a restriction on length, we propose a restriction on {\it shape}. The motivation for our restriction
comes from the fact that if a smooth interpolating curve $F$ has a locally minimal bending energy, then it can
be shown that each segment of $F$, connecting two consecutive interpolation points, is what here is called 
an {\it s-curve}.
In brief, an s-curve is a curve which first turns monotonically in one direction (either
counter-clockwise or clockwise) at most $180^\circ$ and then turns monotonically in the opposite direction at 
most $180^\circ$.
An interpolating curve $F$ is deemed {\it admissible}
if each piece of $F$, connecting two consecutive interpolation points $P_i$ and $P_{i+1}$, is an s-curve. The family
of all admissible interpolating curves is denoted $\Cal A(P_1,P_2,\dots,P_n)$, and we emphasize that our definition
of $\Cal A(P_1,P_2,\dots,P_n)$ includes no restrictions or constraints on length.  Our main result is the following.
\proclaim{Theorem \Tlabel\theoremmain} 
Given any sequence of points  $P_1,P_2,\dots,P_m$ in $\BbR^2$ with $P_i\neq P_{i+1}$, 
the family of admissible interpolating curves
$\Cal A(P_1,P_2,\dots,P_m)$ contains a curve with minimal bending energy.
\endproclaim

An essential sub-problem which arises in the proof of Theorem \theoremmain\ is that of proving the existence of 
an s-curve,
with minimal bending energy, which connects two given unit tangent vectors.  In addition to facilitating our 
proof of Theorem \theoremmain, we anticipate that this sub-problem sits at the core of any numerical
algorithm for solving the general problem,  and with this in mind, we present a thorough analysis of the sub-problem
along with a constructive solution. We mention that in \cite{\GolombJerome} and \cite{\Linner}, the gradient vector 
field approach is employed, but this approach would not apply to $\Cal A(P_1,P_2,\dots,P_m)$ since it is not
an open manifold.


An outline of the sequel is as follows. In section 2, we explain our notation and develop some
basic formulae and properties of rectangular elastica. 
A curve which turns monotonically at most $180^\circ$ in one direction is called a {\it c-curve} and in section
3, we show the existence of an optimal c-curve connecting a unit tangent vector to a line as well as connecting two
unit tangent vectors. Incidentally, line segments are (degenerate) c-curves and c-curves are (degenerate) s-curves. 
In section 4, the uniqueness, or lack thereof, of the optimal c-curves found in section 3 is treated.
The important sub-problem mentioned above, namely the existence of an optimal s-curve connecting two unit tangent vectors, is
primarily solved in section 5, except that one particular case (where the optimal s-curve turns out to be a unique c-curve)
is treated in section 6.  Finally, in section 7, we prove Theorem \theoremmain.

\Section{Notation}
We simplify our notation by using the complex plane $\BbC$ in place of $\BbR^2$.
A {\bf curve}  is a differentiable function 
$f:[a,b]\to\BbC$ whose derivative $f'$ is absolutely continuous and non-zero. The {\bf length} of $f$ is
$\len(f)=\int_a^b |f'(t)|\dd t$. With $L=\len(f)$, let the variables $t\in[a,b]$ and $s\in[0,L]$ be related by
$s=\int_a^t |f'(\tau)|\dd \tau$ and define $F:[0,L]\to\BbC$ by $F(s)=f(t)$. It can be shown that $F$ is
a curve (i.e. $F'$ is absolutely continuous) satisfying $|F'|=1$.  The curve $F$ is called the 
{\bf unit speed curve described by $f$} and is denoted $[f]$. Two curves $f$ and $g$ are said to be {\bf equivalent},
written $f\equiv g$, if $[f]=[g]$.
Since $F'$ is absolutely continuous and $|F'|=1$, it follows that there exists
an absolutely continuous function $\gth:[0,L]\to\BbR$,
unique modulo an additive constant in $2\pi\Bbb Z$, such that $F' = e^{\iMATH \gth}$. 
We refer to $\gth$ at the {\bf direction angle} of $F$, while
the derivative of
$\gth$, denoted $\gk$, is called the {\bf signed curvature} of $F$. Since $\gth$ is absolutely continuous,
it follows that $\gk$ is Lebesgue integrable (see \cite{\Royden, pp. 108--112}). The {\bf turning angle} of $f$, denoted
$\Delta(f)$, is defined by $\Delta(f)=\Delta(F):=\int_0^L \gk(s)\dd s$. Note that the magnitude of the turning angle is
bounded by the $L_1$-norm of $\gk$. If $\gk\geq 0$ (resp. $\gk\leq 0$) almost everywhere in $[0,L]$, then 
$|\Delta(F)|=\norm{\gk}_{L_1}$
and $f$ is called a {\bf left-curve} (resp. {\bf right-curve}). A {\bf c-curve} is a left-curve or a right-curve
whose turning angle has magnitude at most $\pi$. A {\bf u-turn} is a c-curve whose turning angle has magnitude $\pi$.

Given the signed curvature $\gk$ of $F$ and its initial position and direction,
we can recover $F$ as follows:
\newline
Step 1. Define $\gth(s) = \gth_0+\int_0^s \gk(r)\dd r$, $s\in[0,L]$, where $\gth_0 = \arg(F'(0))$.\newline
Step 2. $F(s)=F(0)+\int_0^s e^{\iMATH \gth(r)}\dd r$.\newline
(Here, $\arg$ is defined with the standard range $(-\pi,\pi]$.)
This reconstruction can be used to decide when two curves are close to each other. To see this, suppose $F_1$ is
a unit speed curve having the same length and initial position and direction as $F$. It follows from step 1, that
$\abs{\gth(r)-\gth_1(r)}\leq \norm{\gk-\gk_1}_{L_1}$ and then using the Lipschitz continuity of
the function $r\mapsto e^{\iMATH r}$ in step 2, we obtain
$$
\abs{F(s)-F_1(s)}\leq \int_0^s \abs{e^{\iMATH \gth(r)}-e^{\iMATH \gth_1(r)}}\dd r
\leq \int_0^s \abs{\gth(r)-\gth_1(r)}\dd t
\leq s \norm{\gk-\gk_1}_{L_1}.
\tag\Elabel\closeness$$
Whereas the $L_1$-norm of $\gk$ is necessarily finite, the $L_2$-norm may or may not be finite. When it is finite,
we say that $f$ has finite bending energy, where
the {\bf bending energy} of $f$, denoted $\norm{f}^2$, is essentially the square of the $L_2$-norm of $\gk$:
$$
\norm{f}^2=\norm{F}^2:=\frac14\int_0^L |\gk(s)|^2\dd s.
$$
The constant $\frac14$ has been inserted for later convenience.

A {\bf unit tangent vector} is an ordered pair of complex numbers $u=(u_1,u_2)\in \BbC^2$ such that $|u_2|=1$ and
can be visualized (see Fig. 5.1) as the directed line segment, of unit length, having base-point $u_1$ and direction
$u_2$.
For any $t\in[a,b]$, the {\bf unit tangent vector to} $f$ at $t$, denoted $\T f(t)$, has base-point $f(t)$
and direction $f'(t)/|f'(t)|$, whereby
$$\T f(t):=(f(t),f'(t)/|f'(t)|).$$
The unit tangent vectors $u=\T f(a)$ and $v=\T f(b)$ are called, respectively, 
the {\bf initial} and {\bf terminal} unit tangent vectors of $f$ (see Fig. 5.2),
and we say that $f$ {\bf connects} $u$ to $v$. We also say that $f$ {\bf connects} $u$ to $\gl$ if $\gl$ is the line
through $f(b)$ which is parallel to $f'(b)$ (see Fig. 5.3 where $f_r$ connects $u$ to $\gl$). 
If $g$ is a curve whose initial unit tangent vector equals the terminal
unit tangent vector of $f$, then $[f]$ can be extended by $[g]$ 
obtaining a unit speed curve, denoted $f\splice g$, whose initial and terminal
unit tangent vectors equal those of $f$ and $g$, respectively, and whose bending energy satisfies
$\norm{f\splice g}^2=\norm{f}^2 + \norm{g}^2$.

A {\bf similarity transformation} is a mapping $T:\BbC \to \BbC$ of the form 
$T(z)=c_1 z +c_2$ or $T(z)=c_1 \overline z +c_2$, where
$c_1,c_2$ are complex constants, $c_1\neq 0$. The first form preserves the orientation (left or right) of a curve
while the second form reverses it.
The dilation factor is $\dilpar=|c_1|$, and the effect on a curve $f$ is, as expected, $\len(T\circ f)=\dilpar\, \len(f)$, 
$\norm{T\circ f}^2=\frac1{\dilpar}\norm{f}^2$ and $|\Delta(T\circ f)|=|\Delta(f)|$. If a curve $g$ is equivalent to 
$T\circ f$, then we say that $g$ is {\bf similar} to $f$; in case $\dilpar=1$, $T$ is called a 
{\bf congruency transformation} and we say that $g$ is {\bf congruent} to $f$. Furthermore, we say that $g$ is
{\bf directly} similar (or congruent) to $f$ if $T$ is orientation preserving.

%
%
The curves constructed in this article are formed by line segments (denoted $[A,B]$) and various segments
of rectangular elastica. For the latter, we employ the parameterization $$\Bern(t):= \sin t +\iMATH\x(t),$$ 
where $\x(t)$ is defined by 
$\dsize \frac{d\x}{dt}=\frac{\sin^2 t}{\sqrt{1+\sin^2 t}}$, $\x(0)=0$ (see Figure 6.1a).  
This parameterization is derived simply by substituting $a=1$ and $x=\sin t$ into James Bernoulli's equation 
   $dy = x^2 \dd x / \sqrt{a^4 - x^4}$, and we find it more suitable to the work at hand than the usual unit speed
formulation involving Jacobi's elliptic functions.
Since $\frac{d\x}{dt}$ is even and $\pi$-periodic, it follows that $\x$ is odd and satisfies
$\x(t+\pi)=d+\x(t)$, where 
$$d:=\x(\pi).$$ 
Since the sine function is odd and $2\pi$-periodic, we conclude that
$\Bern(t)$ is odd and satisfies $\Bern(t+2\pi)=\iMATH 2d+\Bern(t)$. We use the notation
$\Bern_{[a,b]}$ to denote the sub-curve $\Bern(t)$, $t\in[a,b]$, and any curve which is similar to $\Bern_{[a,b]}$ is called
a {\bf segment of rectangular elastica}.
For later reference, we mention the following.
$$\align
\abs{\Bern'(t)}=\frac{1}{\sqrt{1+\sin^2 t}}&,\qquad
\frac{\Bern'(t)}{\abs{\Bern'(t)}}=\cos t \sqrt{1+\sin^2 t}+ \iMATH \sin^2 t ,\qquad \gk(t) = 2\sin t,\\
\Delta(\Bern_{[0,b]})&=\int_0^b \gk(t)|\Bern'(t)|\dd t
= 2\cos^{-1}\left[\frac{\cos b}{\sqrt{2}}\right]-\frac\pi2,\\
\norm{\Bern_{[a,b]}}^2 &=\frac14\int_a^b \gk(t)^2|\Bern'(t)|\dd t= \x(b)-\x(a).\\
\endalign
$$
For $t_0\in(0,\pi]$, the segment $\Bern_{[0,t_0]}$ plays an important role in the sequel. In the following lemma,
we establish a connection between the turning angle of $\Bern_{[0,t_0]}$ and the value of $\x(t_0)$.
\proclaim{Lemma \Tlabel\alternate} Let $t_0\in(0,\pi]$ and put $\gth_0=\Delta(\Bern_{[0,t_0]})$. Then
$\dsize \x(t_0)=\frac12\int_0^{\gth_0}\sqrt{\sin\tau}\dd \tau$.
\endproclaim
\demo{Proof} 
Fix $t_0\in(0,\pi]$ and put $\gth=\Delta(\Bern_{[0,t]})=2\cos^{-1}\left[\frac{\cos t}{\sqrt{2}}\right]-\frac\pi2$,
$t\in[0,t_0]$. Then $\frac{d\gth}{dt}=\gk(t)|\Bern'(t)|$, and since $e^{\iMATH\gth}=\Bern'(t)/|\Bern'(t)|$, we have
$\sin\gth=\IM \Bern'(t)/|\Bern'(t)|=\sin^2t$, which implies $\sqrt{\sin\gth}=\sin t$. Hence,
$$
\x(t_0)=\int_0^{t_0}\frac{\sin^2 t}{\sqrt{1+\sin^2 t}}\dd t = \frac12 \int_0^{t_0} \sin(t) \gk(t) |\Bern'(t)|\dd t
=\frac12\int_0^{\gth_0} \sqrt{\sin\tau}\dd \tau.
$$
\qed\enddemo
%
\Section{Existence of optimal c-curves}
%
Given a unit tangent vector $u$ and a line $\gl$, let $C_l(u,\gl)$ denote the set of left c-curves which connect
$u$ to $\gl$. In this section, we consider the problem of finding a curve in $C_l(u,\gl)$ which has minimal bending
energy.
We first consider $C_l(u_0,\gl_d)$, where
$u_0=\T \Bern(0)$ and $\gl_d=\set{z\in\BbC:\IM z = d}$. We will show that $\Bern_{[0,\pi]}$ has
minimal bending energy in $C_l(u_0,\gl_d)$. 
Note that $\norm{\Bern_{[0,\pi]}}^2 = \x(\pi) = d$ and by Lemma \alternate,
we have $d=\x(\pi)=\frac12\int_0^{\pi}\sqrt{\sin\tau}\dd\tau$.

Let $f\in C_l(u_0,\gl_d)$, put $L=\len(f)$, and let $F=[f]$ denote the
unit speed curve described by $f$. Let $\gth$ and $\gk$ be the direction angle and signed curvature
of $F$, respectively,
and note that $\int_0^{ L}\gk(s)\dd s = \pi$ since the turning angle in $C_l(u,\gl_d)$ is $\pi$. Furthermore,
since $F$ originates at $0$ and terminates on $\gl_d$, we have $d=\IM F(L)=\int_0^{ L}\sin\gth\dd s$.
\proclaim{Lemma \Tlabel\smooth} If $\gk$ is continuous and positive, then $\norm{F}^2\geq d$.
\endproclaim
\demo{Proof} 
We adopt the viewpoint that $\gth\in[0,\pi]$ and $s\in[0,L]$ are variables related by
$\gth = \int_0^s \gk(r)\dd r$. The assumptions on $\gk$ ensure that $\gth$ and $s$ are increasing $C^1$ functions
of one another. Noting that $\frac{d\gth}{ds}=\gk(s)$, we observe that
$\norm{F}^2=\frac14\int_0^L \gk(s)^2\dd s = \frac14\int_0^L \gk(s)\frac{d\gth}{ds}\dd s 
=\frac14\int_0^\pi \gk(s)\dd \gth$. Similarly, since $\frac{ds}{d\gth}=1/\gk(s)$, we have
$\int_0^\pi \frac{\sin\gth}{\gk(s)}\dd \gth = \int_0^\pi \sin\gth\frac{ds}{d\gth}\dd \gth
=\int_0^{ L}\sin\gth\dd s=d$.
Now,
$$
d^2 = \left[\frac12\int_0^\pi\sqrt{\sin\gth}\dd \gth\right]^2 
= \left[\int_0^\pi\frac{\sqrt{\sin\gth}}{\sqrt{\gk(s)}}
\frac{\sqrt{\gk(s)}}2\dd \gth\right]^2 
\leq \left(\int_0^\pi\frac{\sin\gth}{\gk(s)}\dd \gth\right)
\left(\frac14\int_0^\pi \gk(s) \dd \gth\right),
$$
by the Cauchy-Schwarz inequality. Hence $d^2\leq d \norm{F}^2$, and therefore $\norm{F}^2\geq d$.
\qed\enddemo
Returning now to the general case, suppose, by way of contradiction, that $\norm{F}^2<d$. Then $\gk$ is a nonnegative
square integrable function satisfying $\int_0^{ L}\gk(s)\dd s = \pi$. It follows that for every $\gep>0$, there 
exists a positive continuous function $\gk_\gep:[0, L]\to[0,\infty)$ such that $\int_0^{ L}\gk_\gep(s)\dd s = \pi$
and $\norm{\gk-\gk_\gep}_{L_2}<\gep$. Let $F_\gep$ be the unit speed curve having signed curvature $\gk_\gep$ and initial
unit tangent vector $u_0$. 
Since $F_\gep$ has turning angle $\pi$, it follows that the terminal unit tangent vector of $F_\gep$ is parallel to
$\gl_d$, but
there is no guarantee that the terminal point $z_\gep=F_\gep( L)$ lies on $\gl_d$. We repair this by
multiplying $F_\gep$ with the positive scalar $c_\gep = d/\IM z_\gep$ obtaining the curve 
$c_\gep F_\gep \in C_l(u_0,\gl_d)$
with $\norm{c_\gep F_\gep}^2 = \frac1{c_\gep}\norm{F_\gep}^2$. 
Since $\norm{\gk-\gk_\gep}_{L_2}\to 0$,
it follows that $\norm{F_\gep}^2\to \norm{F}^2$ as $\gep\to0$.
And since the $L_1$-norm of $\gk-\gk_\gep$ is bounded by a constant multiple of
its $L_2$-norm, it follows from (\closeness) that $c_\gep\to 1$ as $\gep\to 0$. Hence
$\norm{c_\gep F_\gep}^2<d$ when $\gep>0$ is sufficiently small, contradicting Lemma \smooth. Therefore,
$\norm{F}^2\geq d$. We have thus proved the following.
\proclaim{Theorem \Tlabel\Cvl} Let $u_0=\T \Bern(0)$ and $\gl_d=\set{z\in\BbC:\IM z = d}$. Then
$\Bern_{[0,\pi]}$ has minimal bending energy in $C_l(u_0,\gl_d)$.
\endproclaim
\proclaim{Corollary \Tlabel\cCvl} Let $t_1\in (0,\pi)$ and put $u_1=\T \Bern(t_1)$.
Then $\Bern_{[t_1,\pi]}$ has minimal bending energy in $C_l(u_1,\gl_d)$.
\endproclaim
\demo{Proof} If $f\in C_l(u_1,\gl_d)$ has bending energy less than $\Bern_{[t_1,\pi]}$, then
$\Bern_{[0,t_1]}\splice f$ will be a curve in $C_l(u_0,\gl_d)$ with bending energy less than $d$,
contradicting Theorem \Cvl.
\qed\enddemo
\proclaim{Corollary \Tlabel\cCvv} Let $0\leq t_1<t_2\leq \pi$ and put $u_1=\T \Bern(t_1)$, $u_2=\T \Bern(t_2)$. Then
$\Bern_{[t_1,t_2]}$ has minimal bending energy in $C_l(u_1,u_2)$.
\endproclaim
\demo{Proof} If $f\in C_l(u_1,u_2)$ has bending energy less than $\Bern_{[t_1,t_2]}$, then 
$\Bern_{[0,t_1]}\splice f \splice \Bern_{[t_2,\pi]}$ will be a curve in $C_l(u_0,\gl_d)$ with bending energy less than $d$,
contradicting Theorem \Cvl.
\qed\enddemo
\definition{Definition \Tlabel\deflvl} 
Let $\gl$ be a line and $u$ a unit tangent vector whose base-point lies off of $\gl$, and
assume $C_l(u,\gl)$ is nonempty. Let $\delta\in(0,\pi]$ 
be the common turning angle in $C_l(u,\gl)$ and let $t_1\in[0,\pi)$
be such that $\Delta(\Bern_{[t_1,\pi]})=\delta$. There exists a unique similarity transformation $T$ such
that $T\circ \Bern_{[t_1,\pi]}$ belongs to $C_l(u,\gl)$. We define $l(u,\gl)=T\circ \Bern_{[t_1,\pi]}$. In other words,
$l(u,\gl)$ is the unique curve in $C_l(u,\gl)$ which is similar to $\Bern_{[t_1,\pi]}$.
\enddefinition
%
%
\proclaim{Theorem \Tlabel\redhot} 
Let $\gl$ be a line and $u$ a unit tangent vector whose base-point lies off of $\gl$, and
assume $C_l(u,\gl)$ is nonempty. Then $l(u,\gl)$ has minimal bending energy in $C_l(u,\gl)$. Moreover, if
$\delta\in(0,\pi]$ denotes the common turning angle in $C_l(u,\gl)$ and $p$ denotes the orthogonal distance from 
the base-point of $u$ to $\gl$, then 
$$
\norm{l(u,\gl)}^2 =\frac1p\left[\frac12\int_0^\delta \sqrt{\sin\tau}\dd \tau \right]^2.
$$
\endproclaim
\demo{Proof} Let $t_1$ and $T$ be as in Definition \deflvl\ and $u_1$ and $\gl_d$ as in Corollary \cCvl. 
Note that $T$ maps $C_l(u_1,\gl_d)$ onto $C_l(u,\gl)$ and a curve $f$ has minimal bending energy in $C_l(u_1,\gl_d)$
if and only if $T(f)$ has minimal bending energy in $C_l(u,\gl)$. It therefore follows from Corollary \cCvl\ that
$l(u,\gl)=T\circ \Bern_{[t_1,\pi]}$ 
has minimal bending energy in $C_l(u,\gl)$. In order to compute the bending energy of $l(u,\gl)$, recall
that $\norm{l(u,\gl)}^2 = \frac1\dilpar \norm{\Bern_{[t_1,\pi]}}^2$, where $\dilpar$ is the dilation factor in $T$. Since the
orthogonal distance from the base-point of $u_1$ to $\gl_d$ is $d-\x(t_1)$, it follows that $\dilpar=p/(d-\x(t_1))$. And since
$\norm{\Bern_{[t_1,\pi]}}^2=d-\x(t_1)$, we have $\norm{l(u,\gl)}^2 = \frac1p (d-\x(t_1))^2$. By Lemma 
\alternate, $\x(t_1)=\int_0^{\pi-\delta}\sqrt{\sin\tau}\dd \tau$ and hence
$d-\xi(t_1)=\frac12\int_{\pi-\delta}^\pi\sqrt{\sin\tau}\dd \tau=\frac12\int_0^\delta\sqrt{\sin\tau}\dd \tau$.
\qed\enddemo
\remark{Remark \Tlabel\rcurves} The
definitions and results for right c-curves are analogous to those for left c-curves. In brief,
we denote the set of right c-curves connecting $u$ to $\gl$ by $C_r(u,\gl)$, and 
$r(u,\gl)$ is defined the same as $l(u,\gl)$ except that $\delta$ denotes the magnitude
of the common turning angle in $C_r(u,\gl)$ (right curves have a negative turning angle). Theorem \redhot\ then holds
with $C_r(u,\gl)$ and $r(u,\gl)$ in place of $C_l(u,\gl)$ and $l(u,\gl)$, respectively.
\endremark
%
\Section{Uniqueness of optimal c-curves}
%
Having settled the question of existence of an optimal curve in $C_l(u,\gl)$ we now address uniqueness. 
As with existence we start with $C_l(u_0,\gl_d)$, where $u_0$ and $\gl_d$ are as in Theorem \Cvl. 
\proclaim{Theorem \Tlabel\Cvluniq} 
For $i=1,2$, let $F_i:[0,L_i]\to\BbC$ be a unit speed curve in $C_l(u_0,\gl_d)$ 
such that $\norm{F_i}^2=d$ and assume that $F_i$ does not begin or end with a line segment. 
Then $F_1=F_2$.
\endproclaim
Our proof of this
employs the following technical result, which is left as a simple exercise in differential calculus.
\proclaim{Lemma \Tlabel\upside} Let $\nu_1,\nu_2> 0$ and define $H:(0,1)\to(0,\infty)$ by
$H(\mu)=\dfrac{\nu_1^2}\mu + \dfrac{\nu_2^2}{1-\mu}$. 
Then $H$ has a unique minimum at $\mu_0=\nu_1/(\nu_1+\nu_2)$, where $H(\mu_0)=(\nu_1+\nu_2)^2$.
\endproclaim
\demo{Proof of Theorem \Cvluniq} Let $\gth_i$ and $\gk_i$ be the direction angle and signed curvature of $F_i$, respectively. 
Since $F_i$ does not
begin or end with a line segment, we have $0<\gth_i(s)<\pi$ for all $s\in(0,L_i)$, and it follows that $F_i$ can be
reparameterized as $t\mapsto g_i(t)+\iMATH t$, $t\in[0,d]$, where $g_i$ is continuous on $[0,d]$ and 
continuously differentiable on $(0,d)$.
Fix $\gamma\in(0,\pi)$ and let $t\in(0,\pi)$ be such that
$\Delta(\Bern_{[0,t]})=\gamma$. Let $s_i\in(0,L_i)$ be such that $\gth_i(s_i)=\gamma$, and
put $v_i=\T F_i(s_i)$ and $t_i=\IM F_i(s_i)$.  We claim that $t_1=\x(t)=t_2$.
Noting that the turning angle in $C_l(v_i,\gl_d)$ is $\pi-\gamma$ and the orthogonal
distance from the base-point of $v_i$ to $\gl_d$ is $d-t_i$, and
since ${F_i}_{[s_i,L]}$ belongs to $C_l(v_i,\gl_d)$,
we obtain from Theorem \redhot\ and Lemma \alternate\ that $\norm{{F_i}_{[s_i,L_i]}}^2\geq \frac1{d-t_i}(d-\x(t))^2$.
By a similar argument (using right c-curves) we obtain 
$\norm{{F_i}_{[0,s_i]}}^2\geq \frac1{t_i}\x(t)^2$. Therefore,
$$
d=\norm{F_i}^2=\norm{{F_i}_{[0,s_i]}}^2+\norm{{F_i}_{[s_i,L]}}^2\geq \frac1{t_i}\x(t)^2 + \frac1{d-t_i}(d-\x(t))^2.
$$
With $\nu_1=\x(t)$, $\nu_2=d-\x(t)$, $\mu=t_i/d$, and with $H(\mu)$ as in Lemma \upside, we can express
the above inequality as $d\geq \frac1d H(\mu)$, or equivalently, $d^2\geq H(\mu)$. 
By Lemma \upside, $H$ has a unique minimum at
$\mu_0=\x(t)/d$ where $H(\mu_0)=d^2$. But since $d^2\geq H(\mu)$, it must be the case
that $\mu=\mu_0$; therefore $t_i=\x(t)$ as claimed. In terms of the functions $g_1$ and $g_2$, we have
proved that if $g_1'(t_1)=\cot\gamma = g_2'(t_2)$, then $t_1=\x(t)=t_2$. Since, for $i=1,2$, 
$g_i'$ is continuous and decreasing on $(0,d)$, with range $(-\infty,\infty)$, we conclude that $g_1'=g_2'$
on $(0,d)$. Since $g_1(0)=0=g_2(0)$, we have $g_1=g_2$ on $[0,d]$. From this we conclude that $F_1$ and
$F_2$ are equivalent, but since both are unit speed curves, they must be equal.
\qed\enddemo
As an immediate corollary, we have the following.
\proclaim{Corollary \Tlabel\HeyJoe} If $f\in C_l(u_0,\gl_d)$ has minimal bending energy, then $f$ contains a subcurve
which is equivalent to $c+\Bern_{[0,\pi]}$ for some real constant $c\geq 0$.
\endproclaim
Imitating the proof of Corollary \cCvl\ and Theorem \redhot, one easily obtains the following.
\proclaim{Corollary \Tlabel\porcupine} Let $\gl$ be a line and $u$ a unit tangent vector whose base-point lies off of $\gl$, and
assume $C_l(u,\gl)$ is nonempty. Let $\delta\in(0,\pi]$ denote the common turning angle in $C_l(u,\gl)$ and
let $f\in C_l(u,\gl)$ have minimal bending energy. \newline
(i) If $\delta=\pi$, then $f$ contains a subcurve which is congruent to $l(u,\gl)$.\newline
(ii) If $\delta<\pi$, then either $f\equiv l(u,\gl)$ or $f\equiv l(u,\gl)\splice [A,B]$ for some line segment $[A,B]$.
\endproclaim
We have seen in Corollary \cCvv\ that $\Bern_{[t_1,t_2]}$ has minimal bending energy in $C_l(u_1,u_2)$. Using the same
technique as above, one easily obtains the following.
\proclaim{Theorem \Tlabel\cCvvuniq} Let $0\leq t_1<t_2\leq \pi$ and let $u_1=\T \Bern(t_1)$, $u_2=\T \Bern(t_2)$ denote
the initial and terminal unit tangent vectors of the curve $\Bern_{[t_1,t_2]}$, respectively. Let $f$ be a curve
with minimal bending energy in $C_l(u_1,u_2)$. If $f$ is not equivalent to $\Bern_{[t_1,t_2]}$, then
$[t_1,t_2]=[0,\pi]$ and a sub-curve of $f$ is congruent to  $\Bern_{[0,\pi]}$ (i.e. $f$ is equivalent to
$[0,c]\splice (c+\Bern_{[0,\pi]}) \splice [c+\iMATH d,\iMATH d]$ for some real constant $c>0$).
\endproclaim
%
\Section{Optimal s-curves, part I}
%
An {\bf s-curve} is either a c-curve (considered a degenerate s-curve) or
a curve of the form $f=f_1\splice f_2$, where $f_1$ and $f_2$ are c-curves which turn in opposite directions.
Let $u$ and $v$ be two unit tangent vectors and let $S(u,v)$ denote the set of all s-curves which connect
$u$ to $v$. 
In this section and the next, 
we will prove the following.
\proclaim{Theorem \Tlabel\letterman} Let $u$ and $v$ be two unit tangent vectors with distinct base-points. If
$S(u,v)$ is nonempty, then there exists a curve in $S(u,v)$ with minimal bending energy.
\endproclaim
In addition to proving existence, our proof of Theorem \letterman\ will actually describe all optimal curves
in $S(u,v)$. Expecting that the numerical problem of finding an optimal curve in $S(u,v)$ lies at the heart
of future algorithms, we have structured our proof so that it easily translates into a numerical algorithm.

To begin, let $u$ and $v$ be two unit tangent vectors with distinct base-points.
By applying a similarity transformation, if necessary, and possibly a direction reversal (i.e. $S(-v,-u)$ in place
of $S(u,v)$), we
can assume without loss of generality that $u=(0,e^{\iMATH\ga})$ and $v=(1,e^{\iMATH\beta})$, where
$\ga\in[0,\pi]$ and $\abs\gb\leq\ga$ (see Figure 5.1, where $(\ga)$ indicates the direction angle of $u$). 
\newline
\hskip3truecm\epsfysize=2truecm  \epsffile{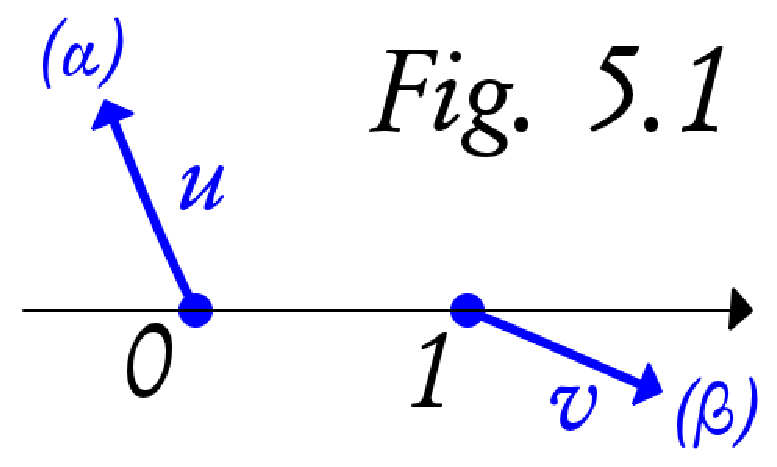}
\newline
We leave it to the reader, as a worthwhile exercise, to verify that in this situation
$S(u,v)$ is nonempty if and only if $\ga<\pi$ and $\gb\geq\ga-\pi$. With that in mind, we proceed assuming that
$\ga\in[0,\pi)$, $\abs\gb\leq\ga$ and $\gb\geq\ga-\pi$.

If $\ga=0$, then $\gb=0$ as well and the line segment $[0,1]$ is the
unique curve (modulo equivalence) in $S(u,v)$ having minimal bending energy. Having dispensed with the trivial case,
we assume henceforth that $\ga>0$. 

Our proof of existence will show that there exists an optimal curve in $S(u,v)$ having one of the following two forms.
\definition{Definition \Tlabel\firstsecond} A curve $f$ is of \newline
(i) {\bf first form} if there exist
$-\pi<t_0<t<\pi$ such that $f$ is directly similar to $\Bern_{[t_0,t]}$,\newline
(ii) {\bf second form} if there exists $c\geq 0$ and $t\in[0,\pi]$ such that $f$ is directly similar to
$$
\Bern_{[-\pi,0]}\splice[0,c]\splice (c+\Bern_{[0,t]}).
$$
\enddefinition
Note that curves of first form do not contain u-turns, while curves of second form do.
While studying a generic right-left s-curve $f\in S(u,v)$,
the following quantities will gradually take on significance, but for easy reference we gather and define them here.
The minimum direction angle $\gga=\min \arg(f')$ is illustrated in Fig. 5.3 for a non-degenerate right-left
s-curve $f$, while $\gga=\gb$ if $f$ is a right c-curve.  The set of all possible angles $\gga$ is denoted $\gG$.
\definition{Definition \Tlabel\Gdefine} For $\gga$ in $\gG:=[\ga-\pi,\gb]\cap (-\infty,0)$, we define
the following:
$$\align
y_1&:=y_1(\gga):=\frac12\int_0^{\ga-\gga}\sqrt{\sin\tau}\dd \tau\qquad(\text{bending energy of }\Bern_{[0,t_1]})\\
y_2&:=y_2(\gga):=\frac12\int_0^{\gb-\gga}\sqrt{\sin\tau}\dd \tau\qquad(\text{bending energy of }\Bern_{[0,t_2]})\\
G(\gga)&:=\frac1{-\sin\gga}(y_1+y_2)^2\qquad\qquad\qquad(\text{lower bound on }\norm{f}^2)\\
\gs(\gga)&:=\cos\gga +\frac{\sin\gga}{y_1+y_2}(\sqrt{\sin(\ga-\gga)}+\sqrt{\sin(\gb-\gga)})\qquad(\text{signed distance})\\
\dilpar(\gga)&:=\frac{-\sin\gga}{y_1+y_2}\qquad\qquad\qquad (\text{dilation factor})
\endalign
$$
Note that, by Lemma \alternate\ (see Fig. 2.1), $y_1$ and $y_2$ can also be expressed as 
$y_1=\xi(t_1)=\norm{\Bern_{[0,t_1]}^2}$ and $y_2=\xi(t_2)=\norm{\Bern_{[0,t_2]}^2}$, where
$t_1,t_2\in[0,\pi]$ are determined by $\gD(\Bern_{[0,t_1]})=\ga-\gga$ and $\gD(\Bern_{[0,t_2]})=\gb-\gga$.
\enddefinition
We mention further that $G(\gga)$ (see Theorem \leftright) is a lower bound on the bending energy of our generic curve $f$,
and $\gs(\gga)$ is a signed distance, which is illustrated in Fig. 5.4. Regarding $\dilpar(\gga)$, we mention that the curves 
$r(u,\gl)$ and $l(\gl,v)$, shown in Fig. 5.4, are similar to $\Bern_{[0,t_1]}$ and $\Bern_{[0,t_2]}$, respectively, with common
dilation factor $\dilpar(\gga)$. The crucial identity relating $G(\gga)$, $\gs(\gga)$ and $\dilpar(\gga)$ is given in Lemma
\Gderiv.

Our constructive proof that $S(u,v)$ contains an optimal curve is broken into three cases which depend on $\ga$ and $\gb$.
To help the reader track these cases, we give here a short description of each case and where in this section
or the next it is treated. 
\definition{Summary \Tlabel\SummarySuv} We assume $\ga\in(0,\pi)$, $\abs\gb\leq\ga$ and $\gb\geq\ga-\pi$.

{\bf Case A:} $\gb=\ga-\pi$.\newline
This case is treated just below and results in an optimal curve of second form.\newline

{\bf Case B:} $\gb\geq 0$ or ($\ga-\pi<\gb<0$ and $\gs(\gb)\geq 0$). \newline
It is shown in Lemma \Gderiv\ that the function $G$ has a 
minimum value $G_{min}$, and in
Corollary \mostofit\ $(vi)$ it is shown that $G_{min}$ equals the minimum bending energy in $S(u,v)$. Each $\gga\in\gG$, 
where $G$ is minimized,
gives rise to an optimal curve in $S(u,v)$, but the form of the optimal curve depends on whether or not $\gga$ is the 
left endpoint of $\gG$.
If $G$ is minimized at the left endpoint $\gga=\ga-\pi$, then it is shown in Corollary \mostofit\ $(iii)$, that the curve $f_{\ga-\pi}$,
which is of second form, is an optimal curve in $S(u,v)$. If $G$ is minimized at any other point $\gga>\ga-\pi$, then it is shown in
Corollary \mostofit\ $(iv),(v)$ that the curve $f_\gga$, which is of first form, is an optimal curve in $S(u,v)$.\newline

{\bf Case C:} $\ga-\pi<\gb<0$ and $\gs(\gb)\leq 0$.\newline
In Theorem \MajorTom, it is shown that the unique curve (modulo
equivalence) in $S(u,v)$ having minimal bending energy is a c-curve of first form.
\enddefinition
\remark{Remark} The reader may note that Case B and Case C have some overlap (namely, when
$\ga-\pi<\gb<0$ and $\gs(\gb)= 0$); this overlap is intentional and serves
as a bridge from Case B to Case C.  The construction under Case C is the better because it yields a unique optimal curve.
\endremark
\demo{Proof of Theorem \letterman\ for Case A}
Assume $\gb=\ga-\pi$. Then $\ga\geq\pi/2$ and $S(u,v)=C_r(u,v)$ (since the common turning angle in $S(u,v)$ is $-\pi$).
Let $\gl$ be the line through $1$ which is parallel to $v$ (and also parallel to $u$) and let 
$P_1$ be the terminal point of $r(u,\gl)$ (see Fig. 5.2).
Since $\ga\geq\pi/2$, it follows that $P_1$ lies on or above the real axis and therefore $f:=r(u,\gl)\splice[P_1,1]$ belongs
to $C_r(u,v)$. Since $S(u,v)=C_r(u,v)\subset C_r(u,\gl)$, it follows from Theorem \redhot\ and Remark \rcurves\ that
$f$ has minimal bending energy in $S(u,v)$. Furthermore, one easily deduces from Corollary \porcupine\ (i) and
Remark \rcurves\ that $f$ is unique modulo equivalence and elongation of u-turns (see remark below).
Note that $f$ is of second form with $t=0$.

Regarding the quantities defined in Definition \Gdefine, the set $\gG$ reduces to the singleton $\gG=\set{\gb}$ 
and since the orthogonal distance from $0$ to $\gl$ is $|P_1-0|=-\sin\gb$, it
follows from Theorem \redhot\ that $\norm{f}^2=d^2/(-\sin\gb)=G(\gb)$. It is easy to see that the distance from 
$P_1$ to $1$ equals $\cos\gb$, and hence $|P_1-1|=\cos\gb=\gs(\gb)$. Lastly, $\dilpar(\gb)=(-\sin\gb)/d$ corresponds
to the dilation factor from $\Bern_{[-\pi,0]}$ to the similar curve $r(u,\gl)$.
\qed\enddemo
%
\epsfysize=5truecm  \epsffile{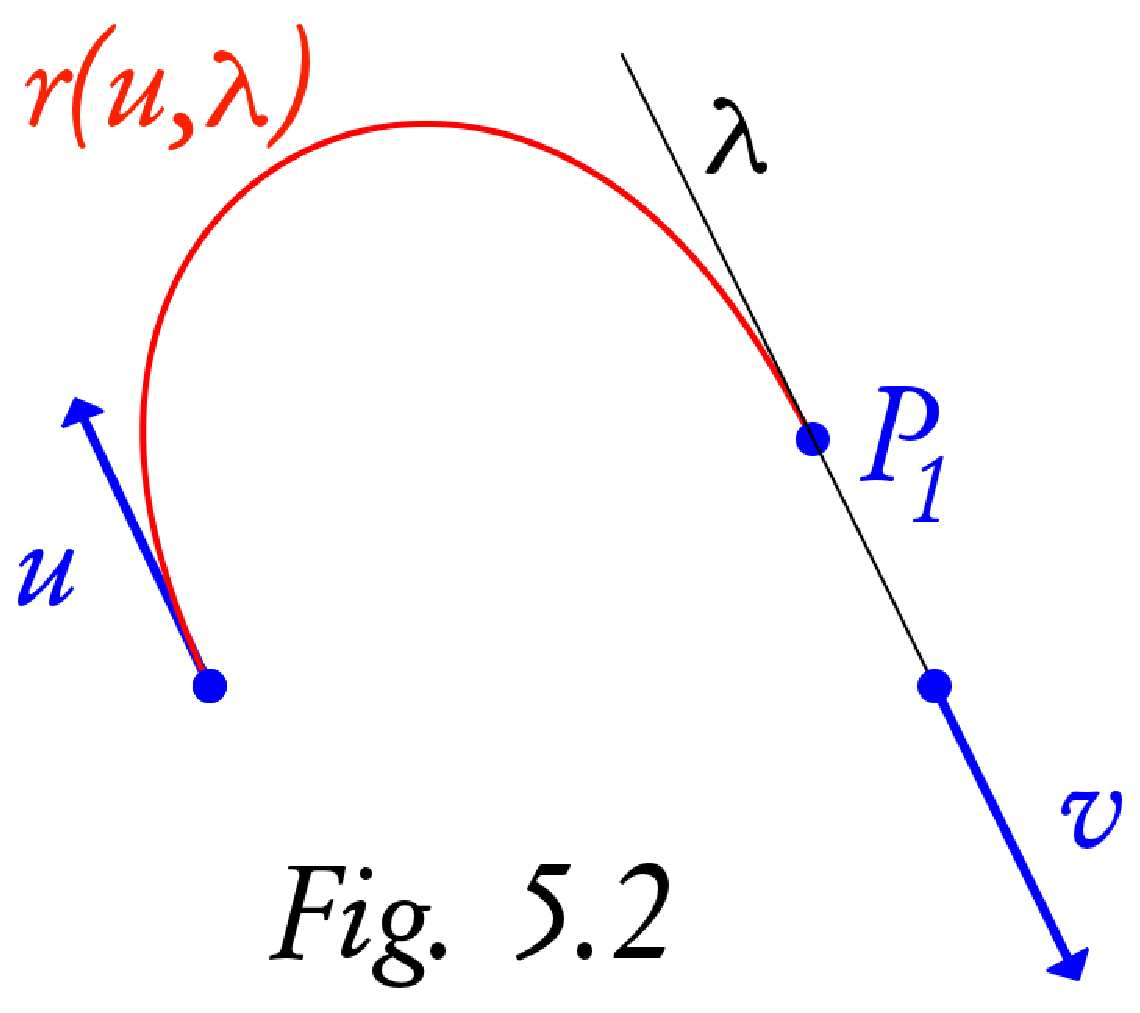}\hskip1truecm
\epsfysize=5truecm  \epsffile{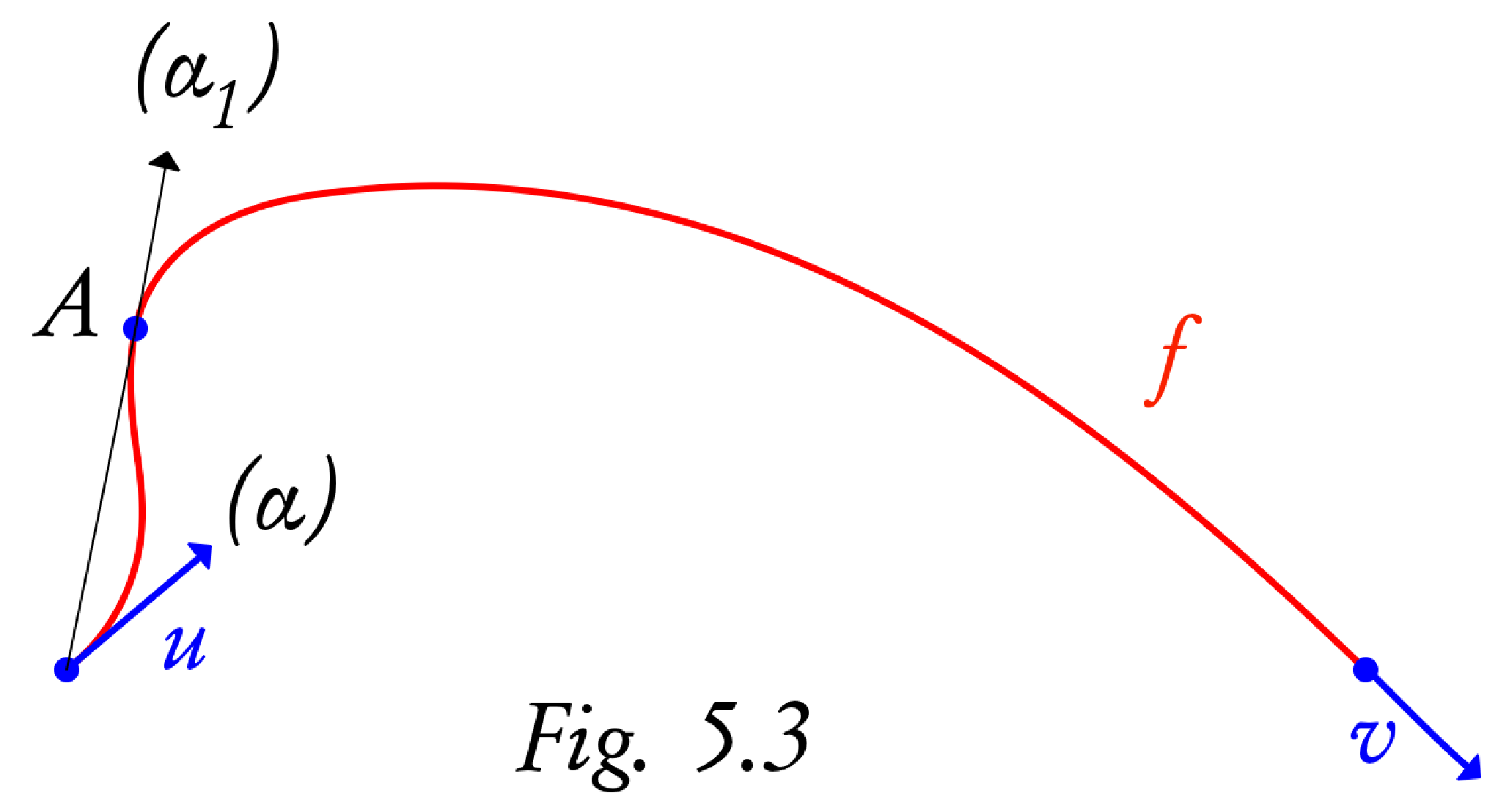}
%
\remark{Remark \Tlabel\Uturn} 
If a curve $f\in S(u,v)$ contains a u-turn (eg. the curve $r(u,\gl)\splice[P_1,1]$ above), then it is 
always possible to
{\it elongate} the u-turn by inserting a pair of congruent line segments before and after the u-turn. 
Although longer, the resulting
curve still belongs to $S(u,v)$ and has the same bending energy as $f$.
\endremark

Having settled the trivial case ($\ga=0$) and Case A, we proceed assuming that
$$
\ga\in(0,\pi),\qquad \abs\gb\leq\ga,\qquad \gb>\ga-\pi.
\tag\Elabel\inforce$$
Our analysis employs an initial partitioning $S(u,v)=S_{lr}'(u,v) \cup S_{rl}(u,v)$, where
$S_{lr}'(u,v)$ (which is nonempty if and only if $\gb<0$) denotes the set of all non-degenerate left-right
s-curves in $S(u,v)$ and $S_{rl}(u,v)$ denotes the set of all right-left s-curves in $S(u,v)$.
\proclaim{Proposition \Tlabel\leftright} If $\gb<0$, then $\norm{f}^2>G(\gb)$ for all $f\in S_{lr}'(u,v)$.
\endproclaim
\demo{Proof} Let $f:[a,b]\to\BbC$ be a non-degenerate left-right s-curve in $S(u,v)$. Set 
$\ga_1=\min_{t\in(a,b]}\arg f(t)$ and let $A=f(s_1)$ be a point where this minimum is attained. It can be shown
that $\arg f'(s_1)=\ga_1>\ga$ and that $f_{[s_1,b]}$ is a right c-curve. Let $u_1$ be the unit tangent vector
$u_1=(0,e^{i\ga_1})$ and let $\gl$ denote the line through $1$ which is parallel to $v$.
Then $g=[0,A]\cup f_{[s_1,b]}$ belongs to $C_r(u_1,\gl)$, and we have
$$\multline
\norm f^2 > \norm{f_{[s_1,b]}}^2=\norm g^2 \geq \norm{r(u_1,\gl)}^2
=\frac1{-\sin\gb}\pr{\frac12 \int_0^{\ga_1-\gb}\sqrt{\sin\tau}\dd \tau}^2\\
>\frac1{-\sin\gb}\pr{\frac12 \int_0^{\ga-\gb}\sqrt{\sin\tau}\dd \tau}^2=\norm{r(u,\gl)}^2=G(\gb).
\endmultline
$$
\qed\enddemo
Although it is not yet apparent, it will eventually be clear that Proposition \leftright\ is all we need
to rule out curves in $S_{lr}'(u,v)$. We now turn our attention to $S_{rl}(u,v)$ and define
subsets $s_\gga^*(u,v)\subset s_\gga(u,v)\subset S_{rl}(u,v)$. For $\gga\in \gG$, let $s_\gga(u,v)$ be the set
of all curves $f\in S_{rl}(u,v)$ whose minimal direction angle, $\min\arg(f')$, equals $\gga$. One easily verifies that
$S_{rl}(u,v)$ partitions as
$
S_{rl}(u,v) = \bigcup_{\gga\in\gG} s_\gga(u,v).
$

If $\gb<0$, then $s_\gb(u,v)$ simplifies to 
$s_\gb(u,v)=C_r(u,v)$ and we also define $s_\gb^*(u,v)=C_r(u,v)$. The definition of $s_\gga^*(u,v)$ is much
more involved when $\gga<\gb$: Let $\gga\in\gG$ with $\gga<\gb$, and let $f\in s_\gga(u,v)$. Then
$f$, being a non-degenerate right-left s-curve, has a well defined inflection line $\gl$ with direction
angle $\gga$ (see Figure 5.4).
\newline
\epsfysize=5truecm  \epsffile{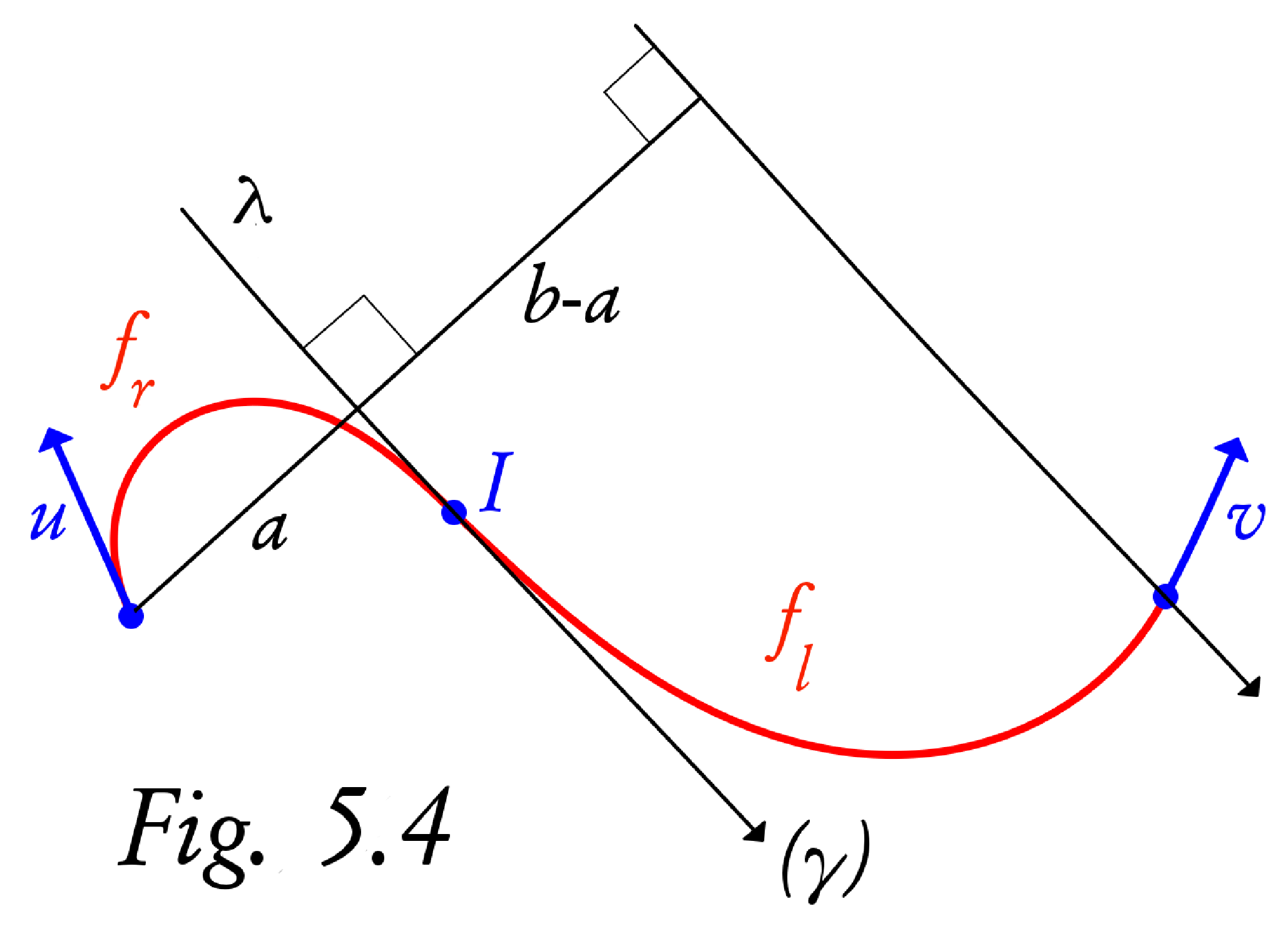}\hskip1truecm
\epsfysize=4truecm  \epsffile{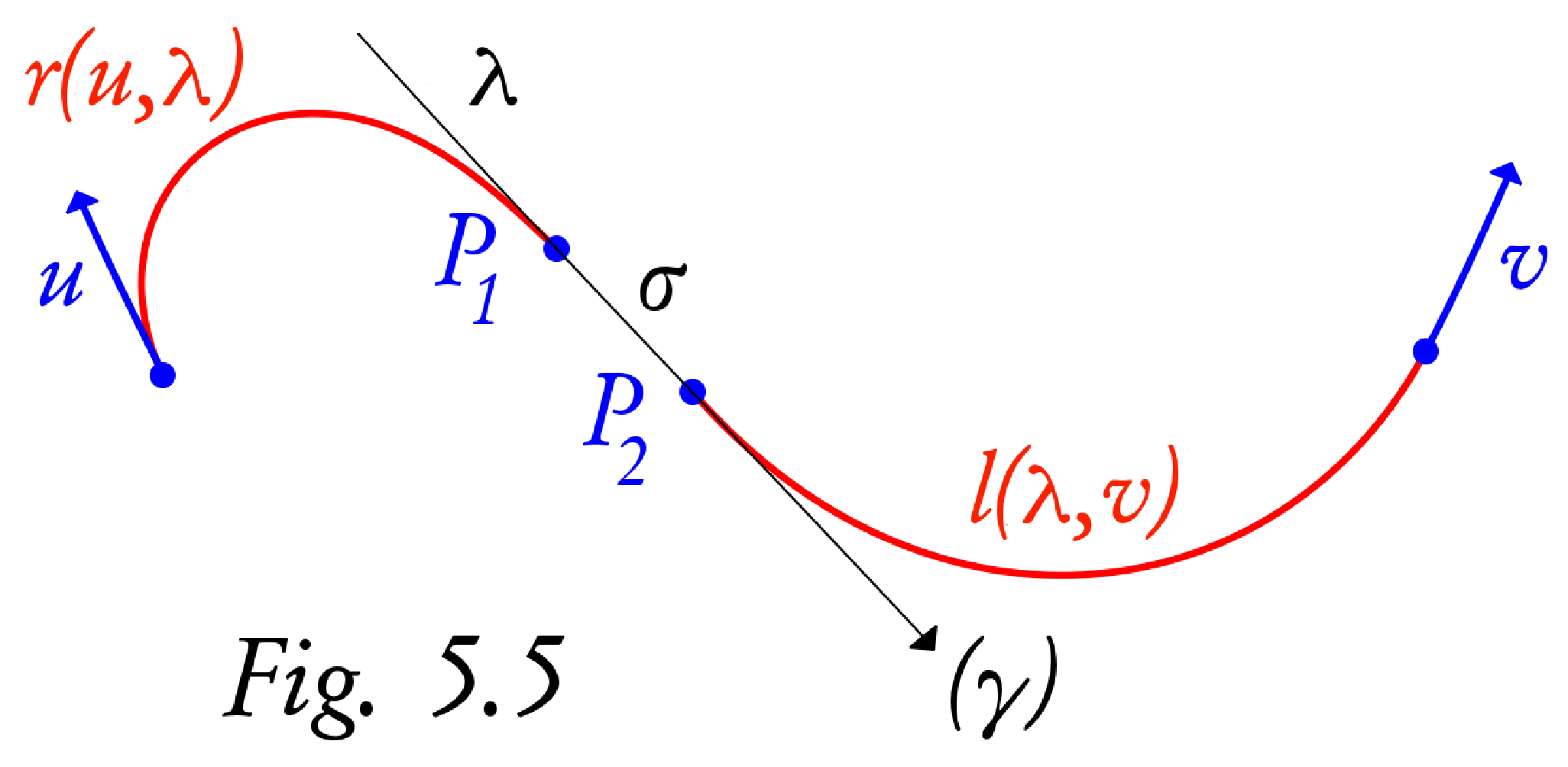}
\newline
Let $a(f)$ denote the orthogonal distance from $0$ to $\gl$. Since $\gl$ necessarily passes between $0$ and $1$, it
follows that $a(f)$ belongs to the interval $(0,b)$, where $b=-\sin\gga$ denotes the orthogonal distance from $0$ to
the line through $1$ and parallel to $\gl$. We define $s_\gga^*(u,v)$ to be the set of all curves $f\in s_\gga(u,v)$
for which $a(f)=-\sin\gga\frac{y_1}{y_1+y_2}$.
\proclaim{Theorem \Tlabel\sgammastar} Assume (\inforce). 
For $\gga\in\gG$ and $f\in s_\gga(u,v)$, the following hold.\newline
(i) $\norm f^2\geq G(\gga)$.\newline
(ii) If $\norm f^2= G(\gga)$, then $f\in s_\gga^*(u,v)$.
\endproclaim
\demo{Proof} We first consider the simpler case when $\gb<0$ and $\gga=\gb$. Since 
$s_\gb^*(u,v)=s_\gb(u,v)=C_r(u,v)$ we need only prove (i). Let $\gl$ be the line through $1$ which is parallel
to $v$. Then $C_r(u,v)\subset C_r(u,\gl)$ and it follows that $\norm f^2 \geq \norm{r(u,\gl)}^2$, since 
$r(u,\gl)$ has minimal bending energy in $C_r(u,\gl)$. Noting that the orthogonal distance from $0$ to $\gl$
is $-\sin\gb$ and the common turning angle in $C_r(u,\gl)$ has magnitude $\ga-\gb$, it follows from
Theorem \redhot\ and Remark \rcurves\ that $\norm{r(u,\gl)}^2=G(\gb)$; hence (i).

Now assume that $\gga\in\gG$ with $\gga<\gb$ and let $f\in s_\gga(u,v)$.
Let the inflection line $\gl$ and the distances $a(f)$ and $b=-\sin\gga$ be as described above (see Fig. 5.4),
and let us write $f=f_r\splice f_l$, where $f_r$ terminates (and
$f_l$ originates) at an inflection point $I$ of $f$. 
Since $f_r\in C_r(u,\gl)$ and $f_l\in C_l(\gl,v)$, it follows
that
$$
\norm f^2 = \norm{f_r}^2+\norm{f_l}^2 \geq \norm{r(u,\gl)}^2+\norm{l(\gl,v)}^2 
= \frac{y_1^2}{a(f)} +\frac{y_2^2}{b-a(f)}.
\tag\Elabel\lateagain$$
Let $H$ be the function defined in Lemma \upside\ with $\nu_1=y_1$, $\nu_2=y_2$ and $\mu=\frac {a(f)}b\in(0,1)$.
Then (\lateagain) can be expressed as $\norm f^2 \geq \frac1b H(\mu)$, and it follows from Lemma \upside\ that
$H(\mu)\geq (y_1+y_2)^2$, with equality if and only if $\mu=y_1/(y_1+y_2)$; that is, if and only if
$a(f)=-\sin\gga\frac{y_1}{y_1+y_2}$. We therefore conclude that $\norm f^2>G(\gga)$ if $f\not\in s_\gga^*(u,v)$ and
$\norm f^2\geq G(\gga)$ if $f\in s_\gga^*(u,v)$, which proves both (i) and (ii). 
\qed\enddemo
For $\gga\in\gG$ with $\gga<\gb$, let $\gl_\gga$ denote the common inflection line in $s_\gga^*(u,v)$, and let
$a_\gga=-\sin\gga\frac{y_1}{y_1+y_2}$ denote the orthogonal distance from $0$ to $\gl_\gga$; furthermore (see Fig. 5.5), 
let $P_1$ and $P_2$ denote the terminal and initial points of $r(u,\gl_\gga)$ and $l(\gl_\gga,v)$, respectively.
If $\gb<0$ and $\gga=\gb$, let $\gl_\gb$ denote the line through $1$ which is parallel to $v$ and let $a_\gga$
and $P_1$ be as defined above (note that $a_\gb=-\sin\gb$ still equals the orthogonal distance from $0$ to
$\gl_\gb$), but for convenience set $P_2 = 1$ (the base-point of $v$). 
We now show that the quantity $\dilpar(\gga)$ corresponds to a common dilation parameter.
\proclaim{Proposition \Tlabel\dilequal} Assume (\inforce) and let $\gga\in\gG$. Then the following hold.\newline
(i) The curve $r(u,\gl_\gga)$ is directly congruent to $\dilpar(\gga) \Bern_{[-t_1,0]}$.\newline
(ii) If $\gga<\gb$, then $l(\gl_\gga,v)$ is directly congruent to $\dilpar(\gga) \Bern_{[0,t_2]}$.
\endproclaim
\demo{Proof} 
Let $T_1$ be the similarity transformation such that $r(u,\gl_\gga)=T_1\circ \Bern_{[-t_1,0]}$ and if $\gga<\gb$,
let $T_2$ be such that $l(\gl_\gga,v)=T_2\circ \Bern_{[0,t_2]}$. 
Since the orthogonal distance from $\Bern(-t_1)$ to the real axis is $\xi(t_1)$
and the orthogonal distance from $0$ to $\gl_\gga$ is $a_\gga$, it follows that the dilation parameter of $T_1$
equals $\frac{a_\gga}{\xi(t_1)}=\frac1{y_1}\frac{y_1}{y_1+y_2}(-\sin\gga)=\dilpar(\gga)$. If $\gga<\gb$, we
see by similar reasoning, that the
dilation parameter of $T_2$ equals $\frac{b-a_\gga}{\xi(t_2)}=\frac1{y_2}(-\sin\gga+y_1\sin\gga/(y_1+y_2))=\dilpar(\gga)$,
where $b=-\sin\gga$ denotes the orthogonal distance from $0$ to the line through $1$ which is parallel to $\gl_\gga$.
\qed\enddemo
The following result shows that the quantity $\gs(\gga)$ corresponds to the signed distance from $P_1$ to $P_2$.
\proclaim{Proposition \Tlabel\sigmameaning} Assume \inforce\ and let $\gga\in\gG$.
Then $\gs(\gga)$ equals the signed distance, in the direction $e^{\iMATH\gga}$, from $P_1$ to $P_2$.
\endproclaim
\demo{Proof} We consider first the case $\gga<\gb$. Let $h$ denote the signed distance in question, and put
$B=P_2-P_1=he^{\iMATH\gga}$ (see Figure 5.5). 
It follows from Proposition \dilequal\ that $f=r(u,\gl_\gga)\splice (l(\gl_\gga,v)-B)$
is directly congruent to $\dilpar(\gga) \Bern_{[-t_1,t_2]}$.
Since the projected distance, in the direction $e^{\iMATH 0}$, from $\Bern(-t_1)$ to $\Bern(t_2)$
equals $\sin t_1 + \sin t_2$, it follows that the projected distance, in the direction $e^{\iMATH\gga}$, from the
initial point to the terminal point of $f$ equals $\dilpar(\gga)(\sin t_1 + \sin t_2)$. Noting that the projected distance,
in the direction $e^{\iMATH\gga}$, from $0$ to $1$ equals $\cos\gga$, we deduce that
$\dilpar(\gga)(\sin t_1 + \sin t_2)=\cos\gga-h$. Solving for $h$ and then employing the identity
$\sin^2\gd = \sin t$, when $\gd=\gD(\Bern_{[0,t]})$ and $t\in[0,\pi]$, yields the desired conclusion $h=\gs(\gga)$.

If $\gb<0$ and $\gga=\gb$, then $t_2=0$ and the above proof, with  $f=r(u,\gl_\gb)$, yields the same conclusion
$h=\gs(\gga)$.
\qed\enddemo
\remark{Remark \Tlabel\twoconsequences} 
Two important consequences of Proposition \dilequal\ and Proposition \sigmameaning\ 
are:\newline
1.  Let $\gga\in\gG$, with $\gga>\ga-\pi$. If $\gs(\gga)=0$, then $f_\gga:=r(u,\gl_\gga)\splice l(\gl_\gga,v)$
has bending energy $G(\gga)$ and
is directly congruent to $\dilpar(\gga) \Bern_{[-t_1,t_2]}$. It follows from the latter that $f_\gga$ is of first form with 
$t_0=-t_1$ and $t=t_2$.\newline
2.  If $\gs(\ga-\pi)\geq 0$, then $f_{\ga-\pi}:=r(u,\gl_\gga)\splice [P_1,P_2]\splice l(\gl_\gga,v)$ 
has bending energy $G(\ga-\pi)$ and is of second form with $c=\gs(\ga-\pi)/\dilpar(\ga-\pi)$ and $t=t_2$.
\endremark
In the following result, we see that $\gs(\gga)$ appears as a factor in the derivative $G'(\gga)$.
\proclaim{Lemma \Tlabel\Gderiv} Assume (\inforce).  The function $G:\gG\to(0,\infty)$ is continuously differentiable,
has a minimum value $G_{min}$, and satisfies $\dsize\frac{d}{d\gga} G(\gga) = \frac1{\dilpar(\gga)^2}\gs(\gga)$ for all
$\gga\in\gG$.
\endproclaim
\demo{Proof} For $\gga\in\gG$, we have
$$\multline
G'(\gga)=\frac{\cos\gga}{\sin^2\gga}(y_1+y_2)^2 - \frac2{\sin\gga}(y_1+y_2)(y_1'(\gga)+y_2'(\gga))\\
=\frac1{\dilpar(\gga)^2}\pr{\cos\gga-\frac{2\sin\gga}{y_1+y_2}(-\frac12\sqrt{\sin(\ga-\gga)}
-\frac12\sqrt{\sin(\gb-\gga)})}=\frac1{\dilpar(\gga)^2}\gs(\gga),
\endmultline
$$
and we note that both $\dilpar$ and $\gs$ are continuous on $\gG$ and $\dilpar$ is positive. If $\gb<0$, then 
$\gG=[\ga-\pi,\gb]$ and it is clear that $G$ has a minimum value. On the other hand, if $\gb\geq 0$, then
$\gG=[\ga-\pi,0)$, but we note that $G(\gga)\to\infty$ as $\gga\to 0^-$; hence $G$ has a minimum value.
\qed\enddemo
In preparation for the main result of this section, we remind the reader that $S(u,v)$ has been partitioned as
$$
S(u,v)=S_{lr}'(u,v) \cup \bigcup_{\gga\in\gG} s_\gga(u,v),
\tag\Elabel\partition$$
where $S_{lr}'(u,v)$ is nonempty only when $\gb<0$.
\proclaim{Corollary \Tlabel\mostofit} Let (\inforce) be in force, and in case $\gb<0$, assume $\gs(\gb)\geq 0$.
The following hold.\newline
(i) If $\gb<0$, then $\norm f^2>G_{min}$ for all $f\in S_{lr}'(u,v)$.\newline
(ii) If $\gga\in\gG$ and $G(\gga)>G_{min}$, then $\norm f^2>G_{min}$ for all $f\in s_\gga(u,v)$.\newline
(iii) If $G(\ga-\pi)=G_{min}$, then $\gs(\ga-\pi)\geq 0$ and the curve $f_{\ga-\pi}$, defined in 
Remark \twoconsequences, is the
unique curve, modulo equivalence and elongation of u-turns, in $s_{\ga-\pi}(u,v)$ with bending energy $G_{min}$.\newline
(iv) Let $\gga\in\gG$ with $\ga-\pi<\gga<\gb$. If $G(\gga)=G_{min}$, then $\gs(\gga)= 0$ and the curve $f_\gga$, defined in 
Remark \twoconsequences, is the
unique curve (modulo equivalence) in $s_\gga(u,v)$ with bending energy $G_{min}$.\newline
(v) If $\gb<0$ and $G(\gb)=G_{min}$, then $\gs(\gb)=0$ and $f_\gb:=r(u,\gl_\gb)$, which is of first form
with $t_0=-t_1$ and $t=0$, is the unique curve
(modulo equivalence) in $s_\gb(u,v)$ with bending energy $G_{min}$.\newline
(vi) The minimum bending energy in $S(u,v)$ is $G_{min}$.
\endproclaim
\demo{Proof} Items (i) and (ii) are immediate consequences of Proposition \leftright\ and Theorem \sgammastar\ (i),
respectively. For (iii), assume
$G(\ga-\pi)=G_{min}$. If $\gs(\ga-\pi)<0$, then it follows from Theorem \Gderiv\ that $G'(\ga-\pi)<0$,
which contradicts the assumption that $G$ attains its minimum at $\ga-\pi$; therefore, $\gs(\ga-\pi)\geq 0$.
It now follows from Remark \twoconsequences\ that $f_{\ga-\pi}$ has bending energy $G_{min}$. Now, suppose
$f\in s_{\ga-\pi}(u,v)$ has bending energy $G(\ga-\pi)=G_{min}$. 
By Theorem \sgammastar\ (ii), $f$ belongs to $s_{\ga-\pi}^*(u,v)$
and writing $f=f_1\cup f_2$, as in the discussion preceding (\lateagain), it follows from (\lateagain) that
$\norm{f_1}^2=\norm{r(u,\gl_{\ga-\pi})}^2$ and $\norm{f_2}^2=\norm{l(\gl_{\ga-\pi},v)}^2$. It can then be deduced
from the results of section 4 
that $f$ is equivalent to $f_{\ga-\pi}$ or can be obtained (equivalently) from
$f_{\ga-\pi}$ by elongation of u-turns. We have thus proved (iii). Turning next to (iv), assume 
$G(\gga)=G_{min}$. Then $G'(\gga)=0$ (since $\gga$ is an interior point of $\gG$) and by Lemma \Gderiv, we have $\gs(\gga)=0$.
It now follows from Remark \twoconsequences\ that $f_\gga$ has bending energy $G_{min}$ and the previous argument
can be applied to show that if $f\in s_\gga(u,v)$ has bending energy $G_{min}$, then $f$ is equivalent to
$f_\gga$ (elongation of u-turns is ruled out since curves in $s_\gga(u,v)$ do not have u-turns). This proves (iv).
For (v), assume $\gb<0$ and $G(\gb)=G_{min}$. If $\gs(\gb)>0$, then $G'(\gb)>0$, by Lemma \Gderiv, which contradicts
the assumption that $G$ is minimized at $\gb$. Therefore, $\gs(\gb)=0$ and it follows that
$r(u,\gl_\gb)$ belongs to $s_\gb(u,v)=C_r(u,v)$ (i.e. $P_1=1$). From Proposition \dilequal\ (i) we have that
$r(u,\gl_\gb)$ is directly congruent to $\dilpar(\gb)\Bern_{[-t_1,0]}$, and from this it is easy to verify that
$\norm{r(u,\gl_\gb)}^2=G(\gb)$. Thus $\norm{r(u,\gl_\gb)}^2=G_{min}$, and we note that
$r(u,\gl_\gb)$ is of first form, with $t_0=-t_1$ and $t=0$. Since the turning angle in $r(u,\gl_\gb)$ has magnitude
less than $\pi$, it easily follows from Theorem \cCvvuniq\ and Remark \rcurves\ that $r(u,\gl_\gb)$ is
the unique curve (modulo equivalence) in $s_\gb(u,v)$ with bending energy $G_{min}$ and the proof of (v) is complete.
We now prove (vi). It follows from (i), Theorem \sgammastar\ and (\partition) that
$\norm f^2\geq G_{min}$ for all $f\in S(u,v)$. Since the function $G$ has a minimum, there exists $\gga\in\gG$ such
that $G(\gga)=G_{min}$, and it then follows from items (iii), (iv) and (v) that $f_\gga$ is a curve in $S(u,v)$ with
bending energy $G_{min}$. This proves (vi).
\qed\enddemo
\noindent As explained in Summary \SummarySuv, Case B of Theorem \letterman\ is a consequence of Corollary \mostofit.
%
\Section{Optimal s-curves, part II}
%
The purpose of this section is to prove the following two results, where we note that Case C (see Summary \SummarySuv)
of Theorem \letterman\ follows from the latter.
\proclaim{Theorem \Tlabel\RocketMan} Let $0\leq t_1<t_2\leq \pi$ satisfy $t_2-t_1<\pi$, and let
$u=\T \Bern(t_1)$ and $v=\T \Bern(t_2)$ be the initial and terminal unit tangent vectors, respectively, to the curve
$\Bern_{[t_1,t_2]}$. Then $\Bern_{[t_1,t_2]}$ is the unique curve (modulo equivalence) in $S(u,v)$ with minimal bending energy.
\endproclaim
\proclaim{Theorem \Tlabel\MajorTom} In the notation of section 5, let $\ga\in(0,\pi)$ and $\gb<0$ satisfy
(\inforce) and suppose  $\gs(\gb)\leq 0$. Then there exist 
$-\pi< t_1<t_2\leq 0$ such that $\Bern_{[t_1,t_2]}$ is directly similar to a curve $f\in S(u,v)$.
Moreover, the curve $f$ is the unique curve (modulo equivalence) in $S(u,v)$ with minimal bending energy.
\endproclaim
For $t\in(0,\pi]$, let $\psi$ and $\phi$, as shown in Figure 6.1b, be the positive angles made by the chord $[0,\Bern(t)]$ and
the segment $\Bern_{[0,t]}$. 
With $\gl$ denoting the tangent line to $\Bern$ at $\Bern(t)$, let $p(t)$ denote the orthogonal  distance from $0$ to $\gl$.
\newline
\epsfysize=5truecm  \epsffile{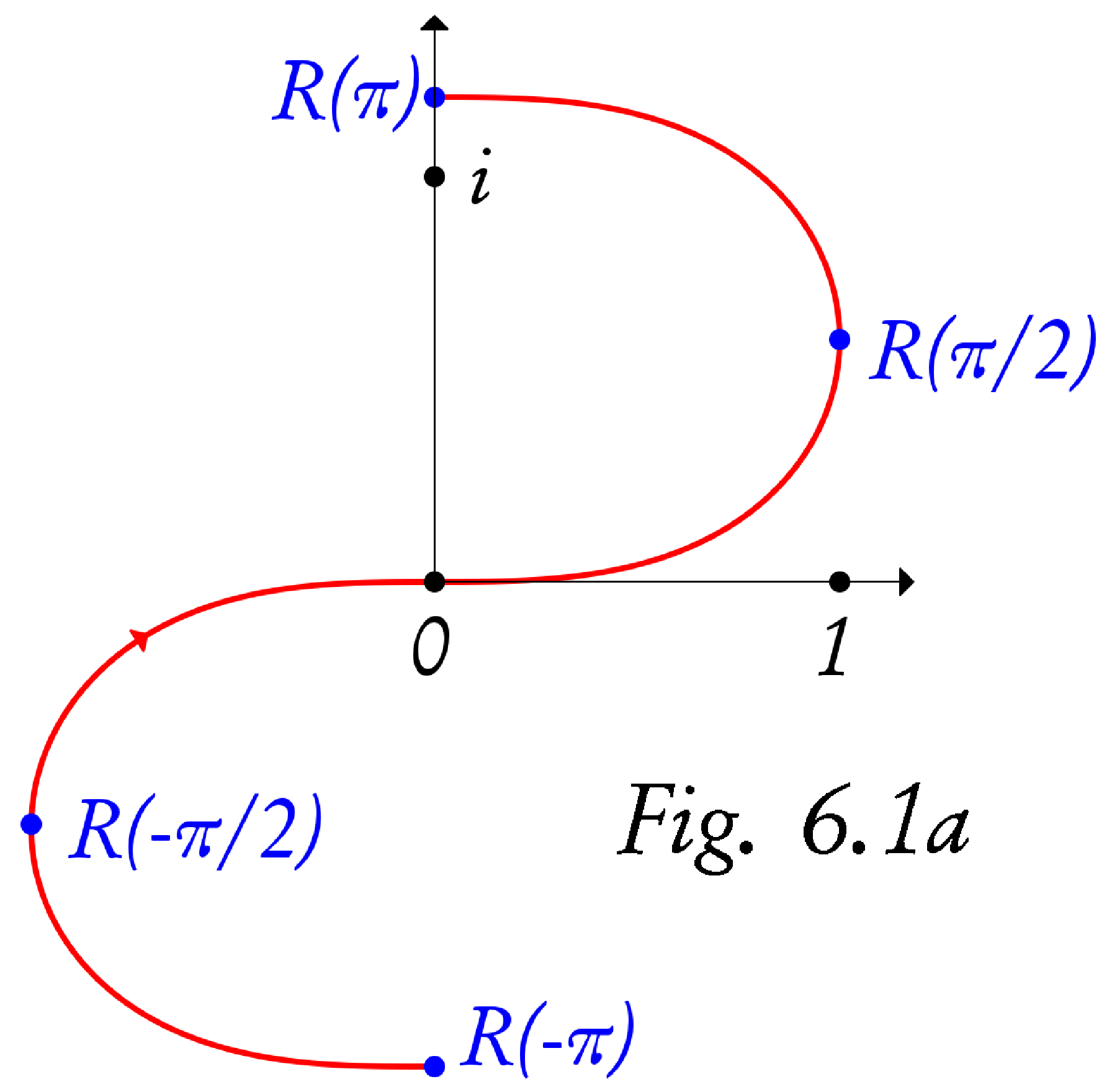}\hskip2truecm
\epsfysize=5truecm  \epsffile{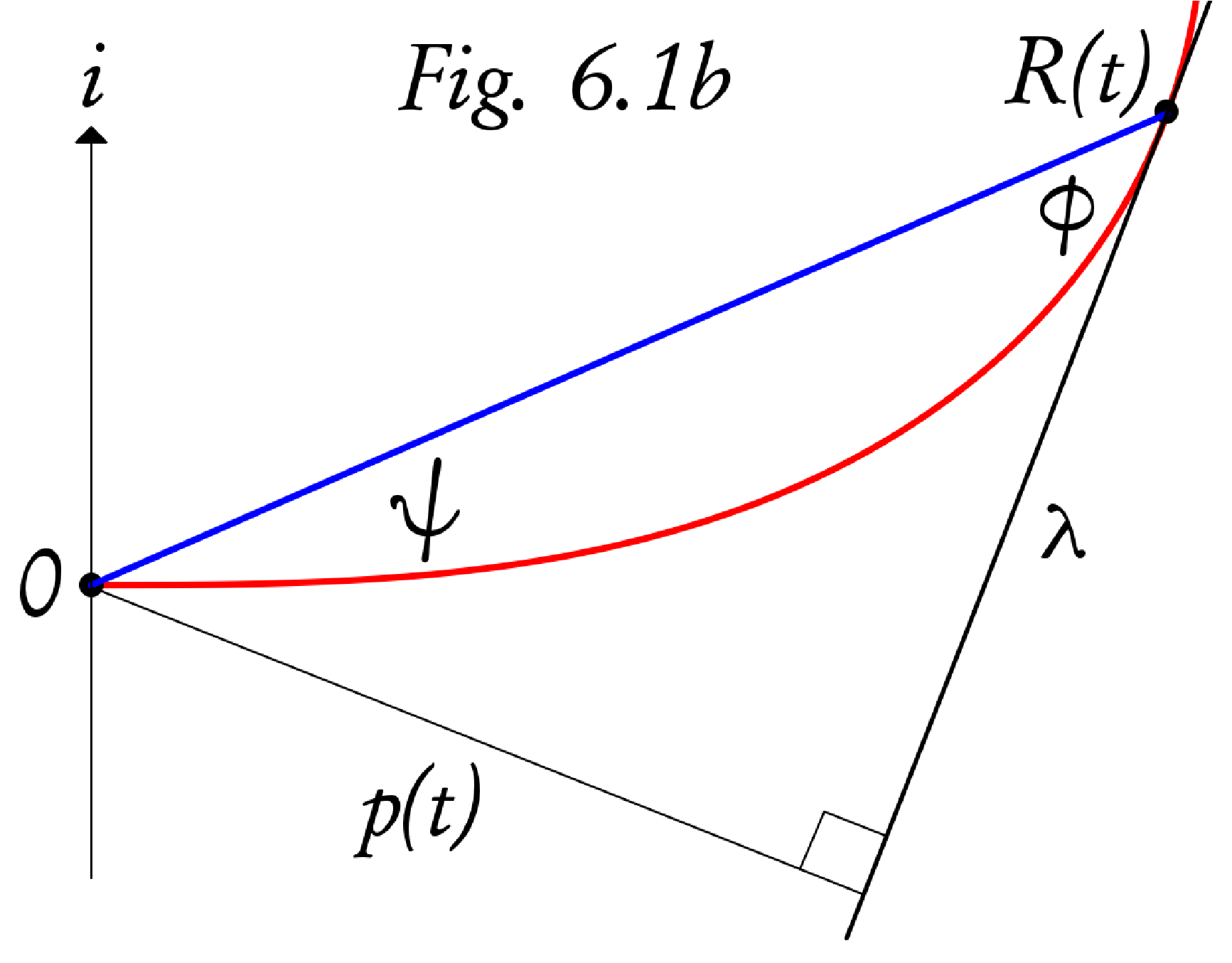}
\newline
\proclaim{Lemma \Tlabel\thbiggerpsi} For $t\in(0,\pi)$, $\phi(t)>\psi(t)$.
\endproclaim
\demo{Proof} We will first show, by way of contradiction, that  $\phi(t)\neq \psi(t)$. Assume $\phi(t)= \psi(t)$ for some
$t\in(0,\pi)$. Let $T(z)=c_1 \overline z + c_2$ be the congruency transformation which interchanges $\Bern(t)$ and $0$, and
set $W=T\circ \Bern_{[t,0]}$, where $\Bern_{[t,0]}$ denotes the reversal of $\Bern_{[0,t]}$.  Since $\phi(t)= \psi(t)$, it follows that $W$ 
belongs to $C_l(\T \Bern(0),\T \Bern(t))$. But
Theorem \cCvvuniq\ asserts that $\Bern_{[0,t]}$ is the unique curve (modulo equivalence) in $C_l(\T \Bern(0),\T \Bern(t))$ with minimal bending energy.
Therefore,
since $W$ and $\Bern_{[0,t]}$ have the same bending energy, they must be equivalent. However, they cannot be equivalent because $W$
begins with a nonzero curvature, namely $2\sin t$, while  $\Bern_{[0,t]}$ begins with curvature $0$. This proves that
$\phi(t)\neq \psi(t)$ for all $t\in(0,\pi)$. While $\phi(\pi)$ and $\psi(\pi)$ both equal $\pi/2$, a simple computation
shows that their derivatives satisfy $-\phi'(\pi)=\psi'(\pi)=1/d$, and it follows that $\phi(t)>\psi(t)$ for $t\in(0,\pi)$
sufficiently close to $\pi$. Since $\phi$ and $\psi$ are continuous, we conclude that  $\phi(t)>\psi(t)$ for all $t\in(0,\pi)$.
\qed\enddemo
\proclaim{Lemma \Tlabel\future} For $t\in[0,\pi)$, $p(t)\xi(t)<(2d-\xi(t))^2$.
\endproclaim
\demo{Proof}  The orthogonal distance $p(t)$ can be formulated as the magnitude of the cross product
$\Bern(t)\times \frac{\Bern'(t)}{|\Bern'(t)|}$ which yields
$$
p(t)=\det\bmatrix  \sin t & \xi(t) \\ 
                   \cos t \sqrt{1+\sin^2 t} & \sin^2t \endbmatrix = \sin^3 t - \xi(t) \cos t \sqrt{1+\sin^2 t},\quad 0\leq t\leq \pi.
$$
We therefore have
$$\multline
p(t)\xi(t)-(2d-\xi(t))^2=(\sin^3 t +4d)\xi(t)-4d^2-\pr{1+\cos t\sqrt{1+\sin^2 t}}\xi(t)^2\\
\leq (\sin^3 t +4d)\xi(t)-4d^2=:g(t),
\endmultline
$$
where the inequality holds since $-1\leq \cos t\sqrt{1+\sin^2 t}\leq 1$.
We note that $g(\pi)=0$ and $\dsize g'(t)=\sin^2t\pr{3\xi(t)\cos t +\frac{\sin^3 t + 4d}{\sqrt{2-\cos t}}}$.
It is clear that $g'(t)>0$ for $t\in (0,\pi/2]$, and for $t\in(\pi/2,\pi)$ (where $-\cos t>0$), we have
$$
g'(t)=(-3\cos t)\sin^2 t\pr{-\xi(t)+\frac{\sin^3 t + 4d}{-3\cos t\sqrt{2-\cos t}}}
\geq(-3\cos t)\sin^2 t\pr{-\xi(t)+\frac{4d}3},
$$
as $0<-\cos t\sqrt{2-\cos t}<1$ on $(\pi/2,\pi)$. Since $0\leq \xi(t)\leq d$, it follows that $g'(t)>0$ for all
$t\in(0,\pi)$ and hence $g$ is increasing on $[0,\pi]$. For $t\in[0,\pi)$, we therefore have
$p(t)\xi(t)-(2d-\xi(t))^2\leq g(t)<g(\pi)=0$, which completes the proof.
\qed
\enddemo
In the following, we again use the notation $S'_{lr}(u,v)$ (resp. $S'_{rl}(u,v)$) for the set of all non-degenerate 
left-right (resp. right-left) s-curves connecting $u$ to $v$.
\proclaim{Lemma  \Tlabel\passormerit} For $t\in(0,\pi)$, the following hold:\newline
(i) If $f\in S_{lr}'(\T \Bern(0),\T \Bern(t))$, then $\norm{f}^2 > \norm{\Bern_{[0,t]}}^2$.\newline
(ii) If $f\in S_{rl}'(\T \Bern(0),\T \Bern(t))$  ends with a left u-turn, then  $\norm{f}^2 > \norm{\Bern_{[0,t]}}^2$.
\endproclaim
\demo{Proof} 
We will employ the notation and results of the previous section, so in order to minimize confusion, we will actually prove
the following equivalent formulations:\newline
$(i')$ If $f\in S_{lr}'(\T \Bern(-t),\T \Bern(0))$, then $\norm{f}^2 > \norm{\Bern_{[-t,0]}}^2$.\newline
$(ii')$ If $f\in S_{rl}'(\T \Bern(-t),\T \Bern(0))$  begins with a right u-turn, then  $\norm{f}^2 > \norm{\Bern_{[-t,0]}}^2$.

Let $T(z)=c_1 z + c_2$ be the similarity transformation determined by $T(\Bern(-t))=0$ and $T(0)=1$, and note that
$T$ brings the configuration $(\T \Bern(-t),\T \Bern(0))$ to the canonical form $(u,v)$ (see Figure 6.2),
where $u=(0,e^{\iMATH\ga})$, $v=(1,e^{\iMATH\gb})$ with  $\ga=\phi(t)$, $\gb=-\psi(t)$. Since $0<\psi(t)<\phi(t)<\pi$,
it follows that (\inforce) holds.
Noting that $r(u,\gl_\gb)=T\circ \Bern_{[-t,0]}$, we see that $\gs(\gb)=0$ and $G(\gb)=\norm{r(u,\gl_\gb)}^2$.
For $(i')$, suppose $f\in S_{lr}'(\T \Bern(-t),\T \Bern(0))$.  Then $T\circ f\in S_{lr}'(u,v)$, and it follows from 
Proposition \leftright\ that  $\norm{T\circ f}^2 > \norm{r(u,\gl_\gb)}^2$. Since $r(u,\gl_\gb)=T\circ \Bern_{[-t,0]}$,
we immediately obtain $(i')$. Now suppose $f\in S_{rl}'(\T \Bern(-t),\T \Bern(0))$ begins with a right u-turn. Then $T\circ f$ belongs to the set
$s_{\ga-\pi}(u,v)$ defined just above Theorem \sgammastar, and it follows from this theorem that
$\norm{T\circ f}^2\geq G(\ga-\pi)$. Since $G(\gb)=\norm{r(u,\gl_\gb)}^2$ and  $r(u,\gl_\gb)=T\circ \Bern_{[-t,0]}$, in order
to establish $(ii')$, it suffices to show that $G(\ga-\pi)>G(\gb)$. From Definition \Gdefine, we have
$G(\gb)=\xi(t)^2/\sin\psi(t)$ and $G(\ga-\pi)=(2d-\xi(t))^2/\sin\phi(t)$. Referring to Figure 6.1b, we see that
$\sin\psi(t) = \xi(t)/|\Bern(t)|$ and $\sin\phi(t)=p(t)/|\Bern(t)|$. Hence
$$
G(\ga-\pi)-G(\gb)
=\frac{|\Bern(t)|}{p(t)}\pr{(2d-\xi(t))^2-p(t)\xi(t)}>0,
$$
by the previous lemma, and this completes the proof of $(ii')$.
\qed\enddemo
\noindent
\epsfysize=5truecm  \epsffile{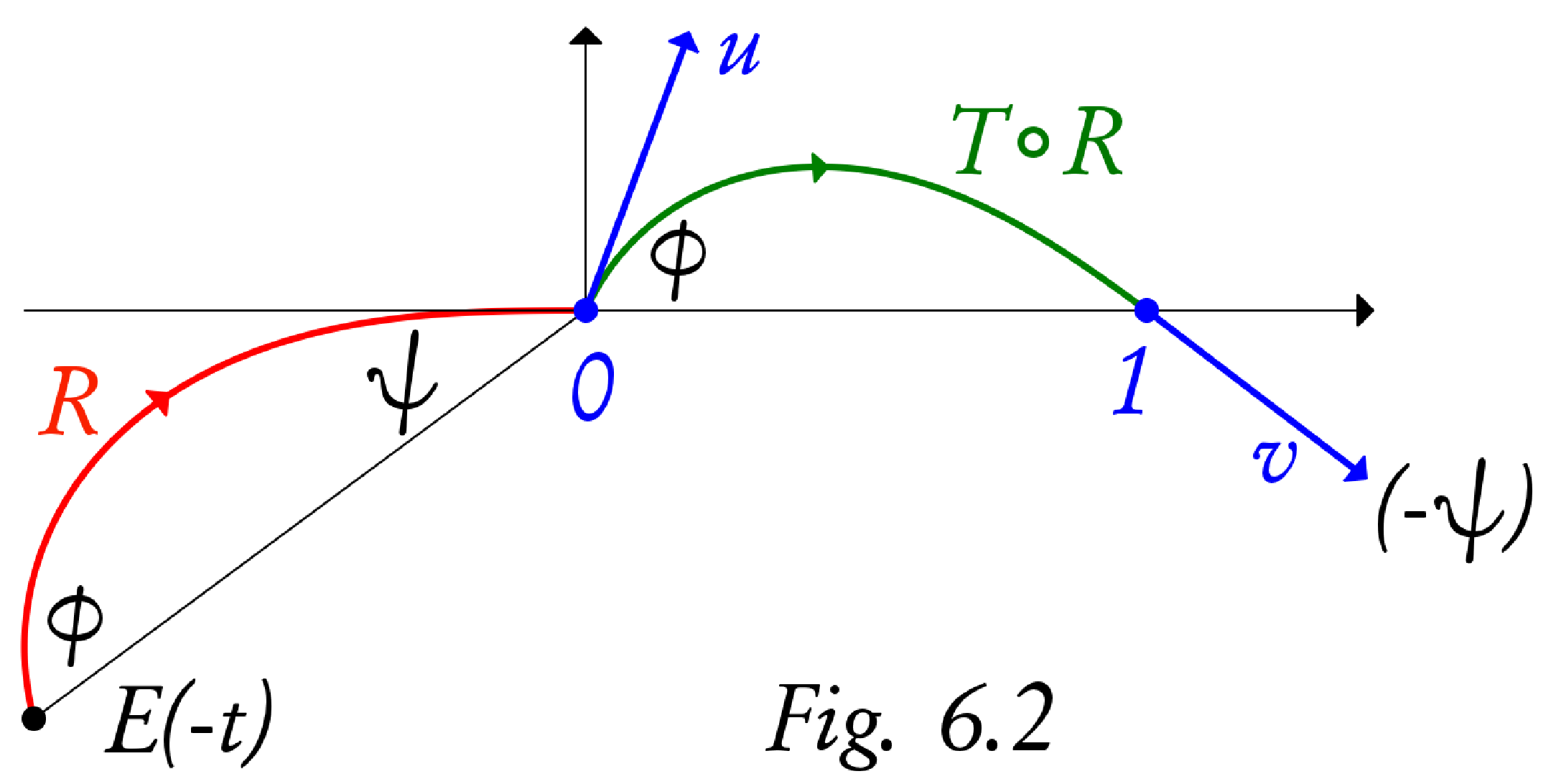}\hskip1truecm
\epsfysize=5truecm  \epsffile{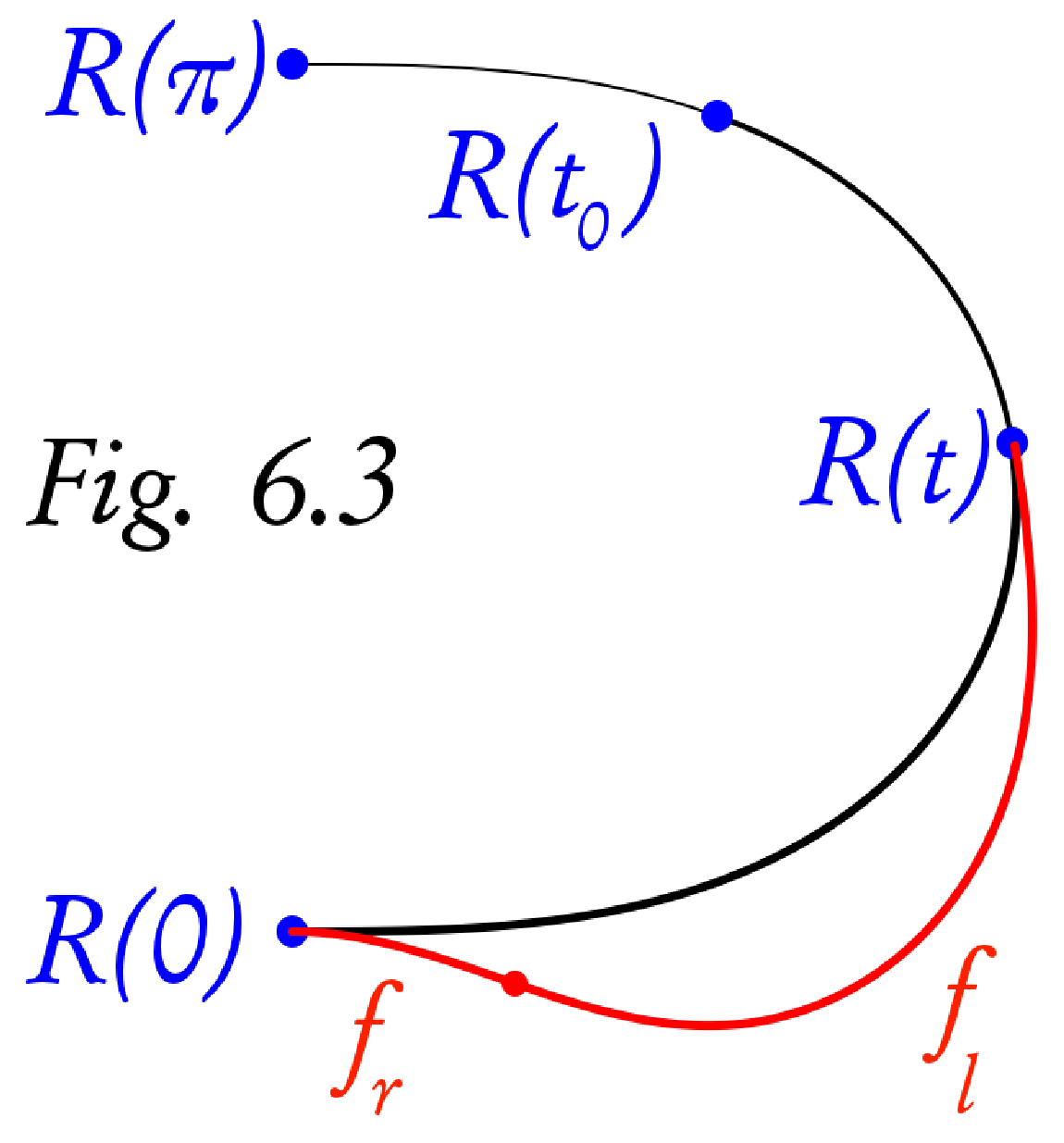}
%
\proclaim{Proposition \Tlabel\halfwaythere} Let $t\in(0,\pi)$. Then
$\Bern_{[0,t]}$ is the unique curve (modulo equivalence) in $S(\T \Bern(0),\T \Bern(t))$ having minimal bending energy.
\endproclaim
\demo{Proof} With Corollary \cCvv, Theorem \cCvvuniq\ 
and Lemma \passormerit\ in view, it suffices to show that $\norm{f}^2 > \norm{\Bern_{[0,t]}}^2$ whenever
$f\in S_{rl}'(\T \Bern(0),\T \Bern(t))$ does not end with a left u-turn. Let $f=f_r\splice f_l$ be as stated, where
$f_r$ is a right c-curve and $f_l$ is a left c-curve (see Figure 6.3). 
Since $f_r$ originates on $\T \Bern(0)$, it follows that there exists
$t_0\in(t,\pi)$ such that $f_l\splice \Bern_{[t,t_0]}$ is a left u-turn. Thus  $f\splice \Bern_{[t,t_0]}$ belongs to
$S_{rl}'(\T \Bern(0),\T \Bern(t_0))$ and ends with a left u-turn. By Lemma \passormerit\ (ii),
$$
\norm{f}^2+\norm{\Bern_{[t,t_0]}}^2=\norm{f\splice \Bern_{[t,t_0]}}^2 > \norm{\Bern_{[0,t_0]}}^2 = \norm{\Bern_{[0,t]}}^2 + \norm{\Bern_{[t,t_0]}}^2,
$$
whence we obtain $\norm{f}^2>\norm{\Bern_{[0,t]}}^2$.
\qed\enddemo
\remark{Remark \Tlabel\remarkhalfway} By symmetry, it follows from Proposition \halfwaythere\ 
that $\Bern_{[t,\pi]}$ is the unique curve (modulo equivalence) in $S(\T \Bern(t),\T \Bern(\pi))$ having minimal bending energy.
\endremark
In the context of the previous section, Proposition \halfwaythere\ asserts the following.
\proclaim{Corollary \Tlabel\bridge} Let direction angles $\ga\in(0,\pi)$, $\gb<0$ satisfy (\inforce) and suppose $\gs(\gb)=0$.
Then $G(\gga)>G(\gb)$ for all $\gga\in[\ga-\pi,\gb)$; that is, $G(\gga)$ is uniquely minimized at $\gga=\gb$.
\endproclaim
\demo{Proof of Theorem \RocketMan} The extreme cases $t_1=0$ and $t_2=\pi$ have been settled in Proposition \halfwaythere\ and
Remark \remarkhalfway, respectively,
so assume  $0< t_1<t_2<\pi$.  By symmetry, and with Corollary \cCvv\ and Theorem \cCvvuniq\  in view, it suffices to show that 
$\norm{f}^2 > \norm{\Bern_{[t_1,t_2]}}^2$ whenever $f$ belongs to $S_{rl}'(\T \Bern(t_1),\T \Bern(t_2))$. Let $f$ be as stated, and let
$\gga\in(-\pi,\pi)$ be the direction angle of $f$ at an inflection point $I$.
\newline
{\bf Case 1:} $\gga\in[0,\pi)$.\newline
Then $f\splice  \Bern_{[t_2,\pi]}$ belongs to $S_{rl}'(\T \Bern(t_1),\T \Bern(\pi))$, and it follows from Remark \remarkhalfway\ that
$\norm{f\splice  \Bern_{[t_2,\pi]}}^2 > \norm{\Bern_{[t_1,\pi]}}^2$, which implies $\norm{f}^2 > \norm{\Bern_{[t_1,t_2]}}^2$.\newline
{\bf Case 2:} $\gga\in(-\pi,0)$.\newline
Since $f$ begins at $\Bern(t_1)$ with a direction angle in $(0,\pi)$, there exists a point $B$ on $f$, between $\Bern(t_1)$
and $I$, where $f$ has direction angle $0$ (see Figure 6.4). 
Let us write $f:=f_1\splice f_2$, where $f_1$ terminates (and $f_2$ originates)
at $B$.
Let $\gl$ be the (horizontal) tangent line to $f$ at $B$, and set $g:=l(\gl,\T \Bern(t_2))$. Since $g$ and
 $\Bern_{[0,t_2]}$ are similar and terminate at the same unit tangent vector, it follows that $g$ originates
at the point of intersection $A$ between $\gl$ and the line segment $[0,\Bern(t_2)]$. Moreover, since $g$ is at a smaller
scale than $\Bern_{[0,t_2]}$, we have $\norm{g}^2>\norm{\Bern_{[0,t_2]}}^2$. Now, it follows from Proposition \halfwaythere\ that
$\norm{[A,B]\splice f_2}^2 > \norm{g}^2$, and therefore
$$
\norm{f}^2>\norm{[A,B]\splice f_2}^2 > \norm{g}^2>\norm{\Bern_{[0,t_2]}}^2>\norm{\Bern_{[t_1,t_2]}}^2.
$$
\qed
\enddemo
\remark{Remark \Tlabel\remarkRocket} By symmetry, Theorem \RocketMan\ remains valid when $-\pi\leq t_1<t_2\leq 0$.
\endremark
\noindent
\epsfysize=6truecm  \epsffile{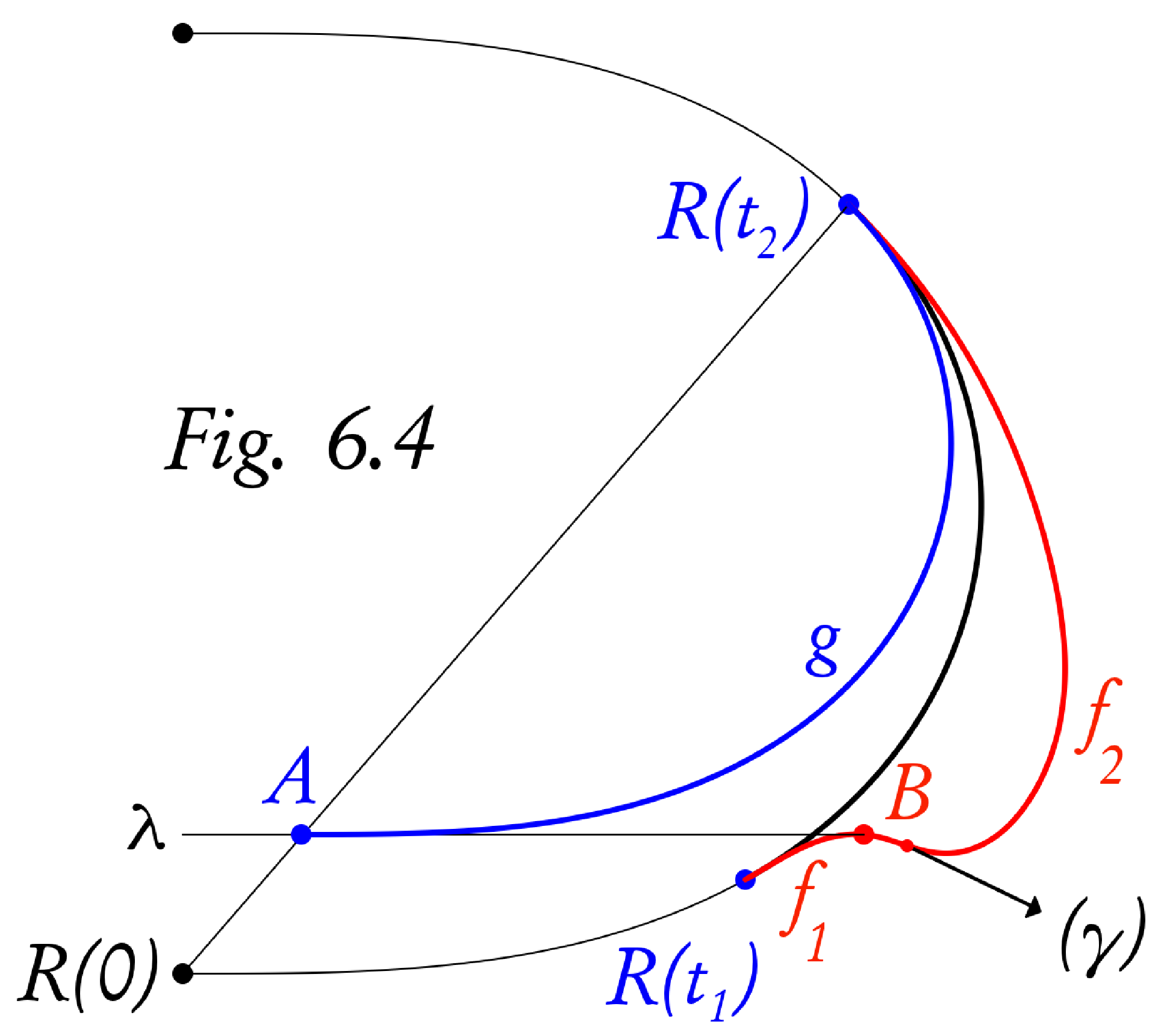}\hskip0.5truecm
\epsfysize=6truecm  \epsffile{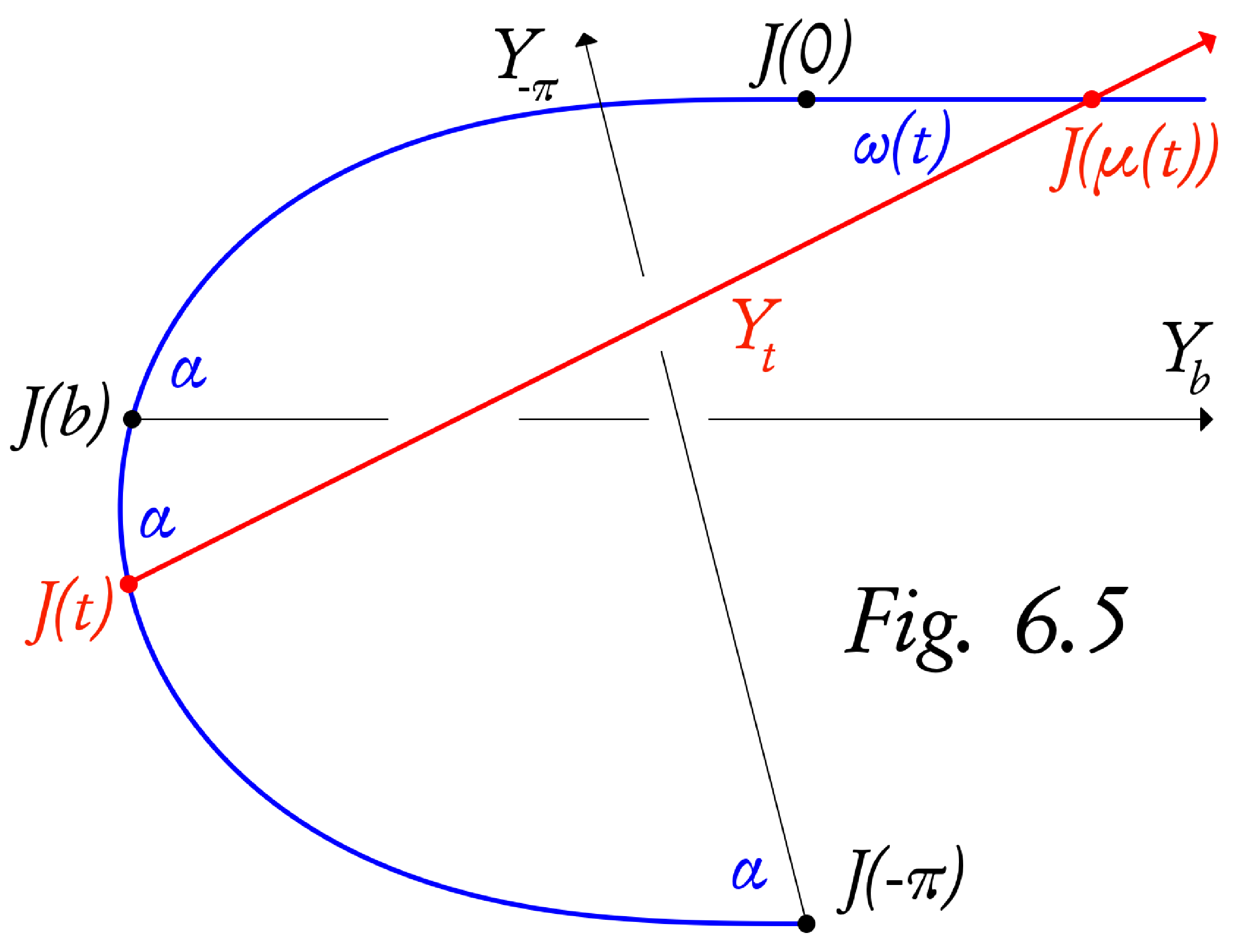}
\proclaim{Lemma \Tlabel\lemmaJ} Define $J:[-\pi,\infty)\to\BbC$ by 
$J(t):=\cases \Bern(t) &\text{if }t\in[-\pi,0]\\
t &\text{if }t>0 \endcases$.
Given positive angles $\ga\geq \gd>0$, with $\ga+\gd<\pi$, there exist $t_1\in(-\pi,0)$ and
$t_2>t_1$ such that the chord $[J(t_1),J(t_2)]$ intersects $J$ with interior angles $\ga$ and $\gd$ at
$J(t_1)$ and $J(t_2)$, respectively.
\endproclaim
\demo{Proof} We refer to Figure 6.5. For $t\in[-\pi,0]$, let $\gth(t)$ denote the direction angle of  $\T \Bern(t)$.
As $t$ ranges from $-\pi$ to $0$, $\gth(t)$  decreases continuously from $\pi$ to $0$,
and it follows that there exists $b\in(-\pi,0)$ such that $\gth(b)=\ga$. For $t\in[-\pi,b)$,
let $Y_t$ denote the ray emanating from $J(t)$ with direction angle $\gth(t)-\ga$ and note that since the
direction angle is positive, $Y_t$ intersects $J$ at a unique point $J(\mu(t))$, where $\mu(t)>t$. Let $\omega(t)$
denote the interior angle, at $J(\mu(t))$,  made when the chord
$[J(t),J(\mu(t))]$ intersects $J$ (the interior angle at $J(t)$ equals $\ga$ by construction). It is clear that
$\omega(t)$ depends continuously on $t\in[-\pi,b)$ and tends to $0$ as $t\to b^-$. We claim that $\omega(-\pi)>\gd$.
If $\mu(-\pi)\geq 0$ (i.e. $J(\mu(-\pi))$ lies on $[0,\infty)$), then $\omega(-\pi)=\pi-\ga$ and the claim follows immediately
from the assumption that  $\ga+\gd<\pi$. On the other hand, if $\mu(-\pi)<0$, then $\ga=\psi(-\mu(-\pi))$ and
$\omega(-\pi)=\phi(-\mu(-\pi))$; hence, by Lemma \thbiggerpsi, $\omega(-\pi)>\ga$ and now the claim follows from the assumption
$\ga\geq\gd$. By the intermediate value property of continuous functions, there exists $t_1\in[-\pi,b)$ such that
$\omega(t_1)=\gd$, and the lemma is proved with $t_2=\mu(t_1)$.
\qed\enddemo
\demo{Proof of Theorem \MajorTom} Put $\gd=-\gb>0$ and note that the hypothesis of Lemma \lemmaJ\ follows from (\inforce), 
and we obtain the conclusion of the lemma. We claim that $t_2\leq 0$. To see this, assume to the contrary that $t_2> 0$.
Let $T(z)=c_1z+c_2$ be the similarity transformation determined by $T(J(t_1))=0$ and $T(J(t_2))=1$. It follows that
$T\circ \Bern_{[t_1,0]} = r(u,\gl_\gb)$ and therefore, by Proposition \sigmameaning, that $\gs(\gb)=\abs{c_1}(t_2-0)> 0$,
which contradicts the assumption that $\gs(\gb)\leq 0$. Therefore, $t_2\leq 0$ and we conclude, from Theorem \RocketMan\ and 
Remark \remarkRocket,
that $f=T\circ \Bern_{[t_1,t_2]}$ is the unique curve (modulo equivalence) in $S(u,v)$ with minimal bending energy.
\qed
\enddemo
%
\Section{Proof of Theorem \theoremmain}
%
For the convenience of the reader let us recall the main theorem
of the paper stated in the introduction (but with $\BbC$ in place of $\BbR^2$).

\proclaim {Theorem {\theoremmain }} Given any sequence of points
$P_1,P_2,\dots,P_m$ in $\BbC$ with $P_i\ne P_{i+1}$, the
family of admissible interpolating curves $\Cal
A(P_1,P_2,\dots,P_m)$ contains a curve with minimal bending
energy.
\endproclaim

The outline of the proof is as follows. First we show that
the family $\Cal A(P_1,P_2,\dots,P_m)$ is non-empty (Prop. 7.1).
With $M$ denoting the infimum of the bending energies of curves in $\Cal A(P_1,P_2,\dots,P_m)$,
it follows that there exists a sequence of curves
$c^1, c^2, c^3, \ldots$ in $ \Cal A(P_1,P_2,\dots,P_m)$
such that $\norm{c^n}^2 \to M$ as $n\to\infty$.

Let $v^n_i$ denote the unit tangent vector
to the curve $c^n$ at the point $P_i$.
Note that the vectors $v_i^1, v_i^2, v_i^3,\ldots$ are all of unit length and
have a common base point $P_i$. Appealing to the Heine-Borel theorem,
and passing to a subsequence if necessary,
we can assume without loss of generality that the sequence
$v_i^1, v_i^2, v_i^3,\ldots$ converges to a unit tangent vector $v_i$,
for $i=1,2,\ldots,m$. To dispel any possible confusion, we mention 
that the sequence of curves $c^1, c^2, c^3,\ldots$
need not converge in any sense--only the unit tangent vectors need converge.

Next we show that for each $i$, the pair $v_i, v_{i+1}$ is
s-feasible (Prop. 7.2), and therefore (by Theorem \letterman), there exists
an s-curve $s_i$  with minimal bending energy in the
family $S(v_i, v_{i+1})$. Joining these pieces together, we construct our candidate
$c=s_1\splice s_2 \splice \cdots \splice s_{m-1}$ which belongs
to $\Cal A(P_1,P_2,\dots,P_m)$.

The proof of Theorem \theoremmain\ is then completed by showing that $\norm{c}^2 = M$.
The proof of this equality uses the fact (Theorem 7.10) that the minimal bending
energy of curves in the family $S(u,v)$ depends continuously
on the directions of $u$ and $v$, and most of the work in 
the current section (Prop. 7.5 -- Prop. 7.9) goes towards establishing this fact.

\proclaim {Proposition \Tlabel \nemures } Under the hypothesis of Theorem \theoremmain,
the family of admissible curves $\Cal A(P_1,P_2,\ldots ,P_m)$ is nonempty.
\endproclaim
\demo{Proof} If we show that there exist
unit tangent vectors $\set{u_j}$, with base-points $\set{P_j}$, such that 
$S(u_j,u_{j+1})$ is nonempty for $j=1,2,\ldots,m-1$, then
$f=f_1\splice f_2\splice \cdots \splice f_{m-1}$, with $f_j\in S(u_j,u_{j+1})$, will be
a curve in $\Cal A(P_1,P_2,\ldots ,P_m)$. We will actually prove a slightly stronger result
in that we will show that $\Cal A(P_1,P_2,\ldots ,P_m,P_1)$ contains a periodic (closed) curve,
where we have tacitly assumed (without loss of generality) that $P_m\neq P_1$. Let $\set{P_j}$
be extended periodically by $P_{j+m}=P_j$.
For $j\in\Bbb Z$, let $w_j=(P_{j+1}-P_j)/|P_{j+1}-P_j|$ be the complex unit in the same direction 
as $P_{j+1}-P_j$ (see Figure 7.1). We then set $u_j=(P_j,z_j)$, where $z_j$, the direction of $u_j$, is defined as follows:\newline
If $w_{j-1}+w_j$ is nonzero, then $z_j$ is the complex unit in the same direction
as $w_{j-1}+w_j$; otherwise, $z_j=e^{i\pi/2} w_j$. Since $\set{u_j}$ has inherited the periodicity of $\set{P_j}$
(namely $u_{j+m}=u_j$), in order to complete the proof, it suffices to show that $S(u_j,u_{j+1})$ is nonempty for all 
$j\in\Bbb Z$. 
\newline
\hskip2cm\epsfysize=4truecm
\epsffile{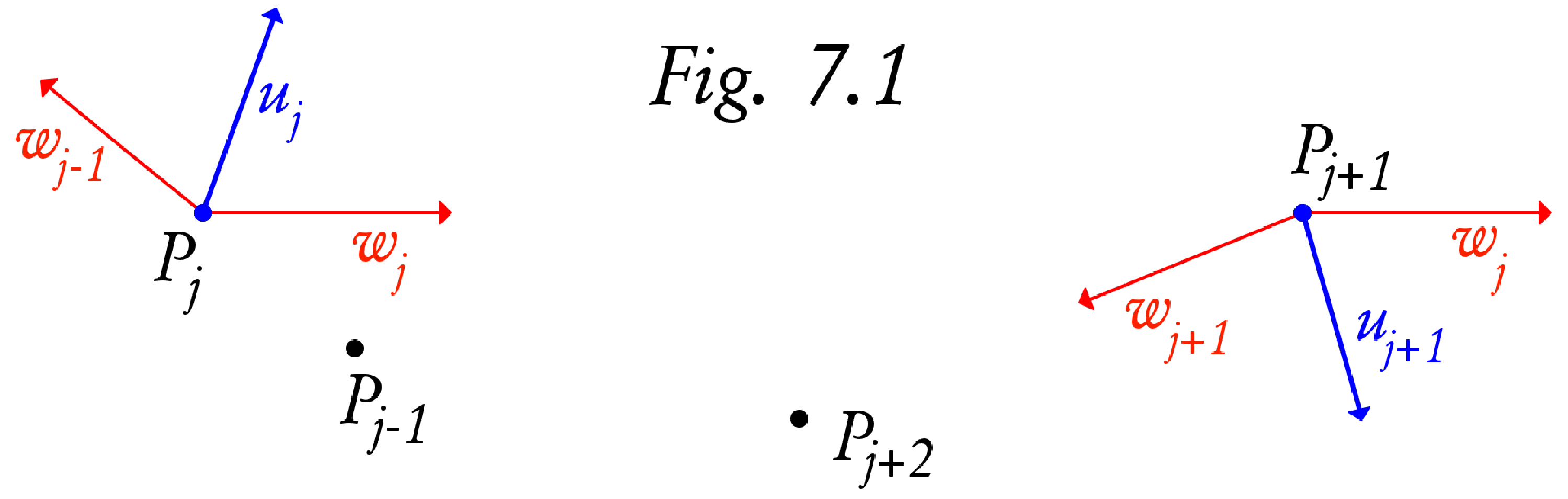}\hskip0.5truecm
\newline
Fix $j\in\Bbb Z$. By applying a rigid motion, if necessary, we can assume without loss of generality that 
$P_j=0$ and that $P_{j+1}$ lies on the positive real axis. That
$S(u_j,u_{j+1})$ is nonempty can now be established simply by showing that both directions $z_j$
and $z_{j+1}$ have nonnegative real part. Note that $w_j=1$. If $w_{j-1}\neq -1$, then 
$z_j$ is the complex unit in the same direction as $w_{j-1}+1$ and hence
$\RE z_j\geq 0$; otherwise, $\RE z_j= \RE e^{i\pi/2}=0$. By the same reasoning,
one sees that $\RE z_{j+1}\geq 0$.
\qed\enddemo

Let $u,v$ be s-feasible ($S(u,v)\ne \emptyset $) unit tangent
vectors. Denote by $\Energy(u,v)$  the minimum of the bending
energy in $S(u,v)$. By the previous sections  (Theorem
{\letterman}) this minimum is assumed by a curve in $S(u,v)$.

We will need the following proposition which shows that the limit
of s-feasible vectors is also s-feasible assuming that the bending
energy is bounded.

\proclaim {Proposition \Tlabel \limitspline } Let $P_u\ne P_v$ be
different points of $\BbC$. Let $u_n, v_n$ be s-feasible unit
tangent vectors with base-points $P_u$ and $P_v$ respectively,
such that $\lim u_n=u$ and $\lim v_n=v$. If $\set{\Energy(u_n,v_n)}$
is bounded then $u,v$ is also s-feasible.
\endproclaim

\demo {Proof} Assume $\set{\Energy(u_n,v_n)}$ is bounded.
Without  loss of generality we can assume that
$P_u=0\in \BbC$ and $P_v=1\in \BbC$. Let $\alpha ,\alpha _n,\beta,\beta_n\in(-\pi,\pi]$
be the direction angles of $u,u_n,v,v_n$, respectively (see Figure 5.1). Since 
$u_n,v_n$ are s-feasible configurations, we must have $\alpha _n,\beta_n\in(-\pi,\pi)$.
Moreover, it is easy to see that if
with $|\alpha _n|\to \pi$ or $|\beta_n|\to\pi$,  then $\Energy(u_n,v_n) \to \infty
$; therefore $\alpha,\beta\in(-\pi,\pi)$. 

If $\alpha =0$ then $u,v$ is s-feasible, regardless of $ \beta \in (-\pi, \pi )$.
So assume $\alpha \ne 0$. Without loss of generality, we may assume that
$\alpha >0$ and $\alpha _n>0$ for all $n\in \Bbb N$. Similar to the exercise assigned
to the reader at the beginning of section 5, we leave it to the reader to verify
that $u_n, v_n$ is s-feasible if and only if  $\beta_n\in [\alpha_n-\pi, \pi)$. Since $\beta
\in (-\pi, \pi )$, it follows that $\beta \in [\alpha -\pi, \pi)$, and therefore the limit 
configuration $u,v$ is s-feasible. 
\qed\enddemo

Next, we  show that the bending energy $\Energy(u,v)$ is
continuous in $u$ and $v$. For this we will need some preparation.

Let $u=(0,e^{i\ga})$ and $v=(1,e^{i\gb})$ be two unit tangent vectors with direction angles
$\ga,\gb\in(-\pi,\pi]$ and for comparison, let $\bar u=(0,e^{i\bar\ga})$ and 
$\bar v=(1,e^{i\bar\gb})$ be unit tangent vectors
with the same base-points as $u$ and $v$, but possibly different directions. The {\bf diameter}
of a curve $f:[a,b]\to\BbC$ is defined by $\diam(f) = \max_{t,\tau} |f(t)-f(\tau)|$.
\proclaim{Lemma \Tlabel\ccurveone} With the notations introduced above,
let $D>0$ and $\eta\in(0,\pi/4)$ be given and assume that 
$\ga,-\gb\in(\eta,\pi-\eta)$ are such that there exists a curve
$f:[a,b]\to\BbC$ in
$C_r(u,v)$ with $\diam(f)\leq D$.  Then, for every
$\gep >0$ there is a $\gd_1=\delta_1 (\gep, D,\eta )>0$ (depending
only on $\gep, D,\eta $) such that if $|\bar \alpha -\alpha |<\delta_1$ and
$C_r(\bar u, v)$ is nonempty, then there is a curve $c\in C_r(\bar u,v)$ 
such that $||c||^2 \leq ||f||^2+ \gep$ and $\diam(c)\leq 2D+1$.
\endproclaim
\demo{Proof} 
We will describe, in two cases, how to modify the curve $f$ near the base point of $u$ to obtain a new 
curve $c$ in $C_r(\bar u, v)$ satisfying $\norm{c}^2\leq \norm{f}^2 +\gep$.

\noindent{\bf Case 1:} $0<\bar \alpha <\alpha $.
\epsfysize=4truecm
\epsffile{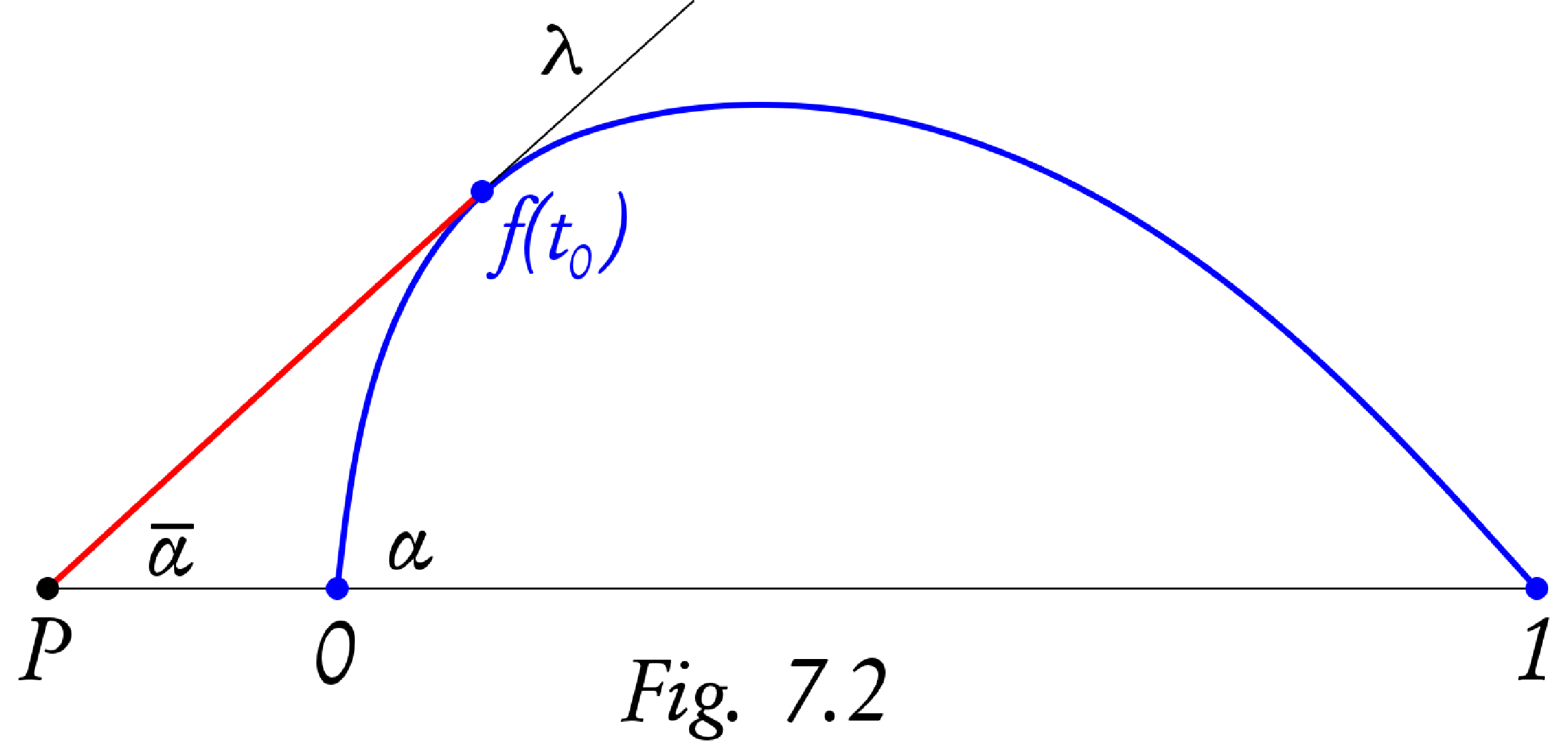}\hskip0.5truecm
\newline
Draw a line with direction angle $\bar \alpha $ through a distant point
on the negative real axis, and
then translate it horizontally towards $f $ until it
makes first contact, obtaining a line $\gl$ (see Figure 7.2). Let $P$ denote the intersection point 
of $\gl$ with the negative real axis and let $f(t_0)$ be a point on $\gl$.
Set $ c_1=[P,f(t_0)]\splice f_{[t_0,b]}$ and note that $c_1$ is similar to a curve $c$ in $C_r(\bar u,v)$
having bending energy $\norm{c}^2=(1-P)\norm{c_1}^2$, while
$\norm{c_1}^2 =\norm{ f_{[t_0,b]}}^2 \leq \norm{f}^2$.
It is easy to see that there is a $\delta_1=\gd_1(\gep,\eta , D)\in(0,\eta/2)$
such that if $\alpha -\delta _1<\bar \alpha <\alpha$, then
the distance from $P$ to $0$ is less than the minimum of $1$ and
$\frac{\gep}{||f||^2}$; hence
$||c||^2\leq (1-P)\norm{f}^2 \leq ||f||^2+\gep$ and $\diam(c)\leq 2D+1$.\newline
\noindent{\bf Case 2:}  $\bar \alpha >\alpha $.
\epsfysize=4truecm
\epsffile{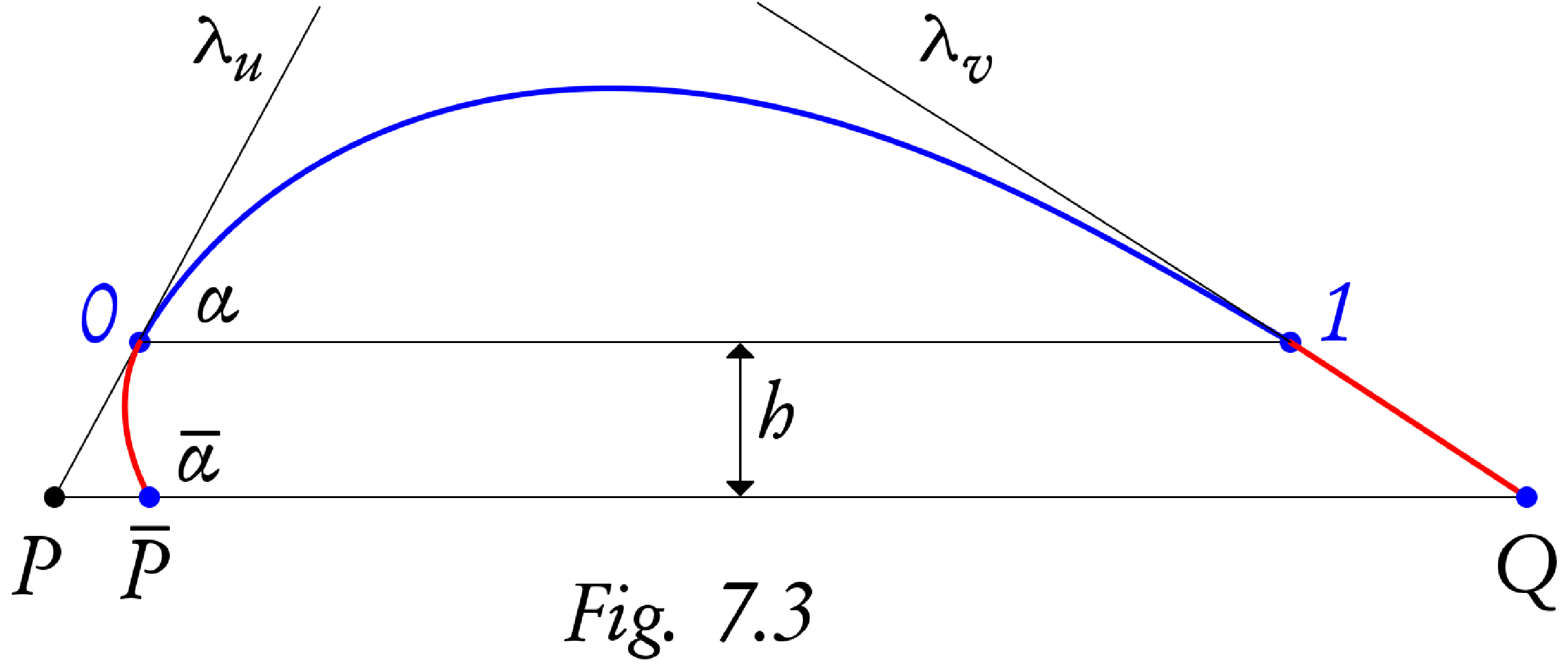}
\newline
We obtain our curve $c$ in two steps (see Figure 7.3).
First, draw the  tangent lines $\gl_u$ and $\gl_v$ to $f$ at $0$
and $1$, respectively.  For a given $h>0$, let $P$ and $Q$ be the points on
$\gl_u$ and $\gl_v$, respectively, having $\IM P = \IM Q = -h$. Since, $\alpha ,-\beta \in(\eta,\pi-\eta) $,
there exists $h=h(\eta,\gep)\in(0,1/2)$ such that the resulting points $P$ and $Q$ will satisfy $|P|<1/9$,
$|Q-1|<1/9$ and $Q-P<\frac{||f||^2+\gep}{||f||^2+\gep/2}$.
Let $h$, $P$ and $Q$ be thus fixed.

Second, replace the line segment $[P,0]$ by an arc of a circle
which emanates from a point $\bar P$ on $[P,Q]$, with direction angle $\bar\ga$,
and terminates at $0$ with direction angle $\ga$, as shown in Figure 7.3.
If $\bar\ga$ is sufficiently close to $\ga$, one can construct such an arc
easily and we leave the details to the reader. With $g$ denoting this circular
arc, we define $c_1=g\splice f \splice [1,Q]$, which is similar to
a curve $c$ in $C_r(\bar u,v)$, and we note that
$\norm c^2 = (Q-\bar P)\norm{c_1}^2$, while $\norm{c_1}^2 = \norm{f}^2 + \norm{g}^2$.
The radius of the circle containing $g$ tends to $\infty$ as $\bar \ga$ tends to $\ga$, and therefore
$\norm g^2$ can be made arbitrarily small by taking $\bar\ga$ sufficiently close to $\ga$.
Consequently, there exists $\gd_1=\gd_1(\eta,\gep)\in(0,\eta/2)$ such that
if $\alpha <\bar \alpha <\alpha +\delta _1$, then $\norm{g}^2<\gep/2$ and consequently we have
$\diam(c)<2D+1$ and $\norm c^2  < (Q-P)(\norm f^2 + \norm g^2) \leq \norm f^2 +\gep$.
\qed\enddemo
\remark{Remark \Tlabel\ccurvetwo} Under the same hypothesis as Lemma \ccurveone, we can prove in a similar
manner that for every
$\gep >0$ there is a $\gd_2=\delta_2 (\gep, D,\eta )>0$ (depending
only on $\gep, D,\eta $) such that if $|\bar \gb -\gb |<\delta_1$ and
$C_r(u, \bar v)$ is nonempty, then there is a curve $c\in C_r(u, \bar v)$ 
such that $||c||^2 \leq ||f||^2+ \gep$ and $\diam(c)\leq 2D+1$.
\endremark

\proclaim{Proposition \Tlabel\ccurve } Under the same hypothesis as Lemma \ccurveone,
for every $\gep >0$ there is a $\gd_3=\delta_3 (\gep, D,\eta )>0$ (depending
only on $\gep, D,\eta $) such that if 
$|\bar \alpha -\alpha |,|\bar \gb -\gb |<\delta_3$ and
$C_r(\bar u, \bar v)$ is nonempty, then there is a curve $c\in C_r(\bar u,\bar v)$ 
such that $||c||^2 \leq ||f||^2 + \gep$.
\endproclaim
\demo{Proof} Let $\gep>0$ be given, and let $\gd_3$ be the minimum of $\eta/2$,
$\gd_1(\gep/2,D,\eta)$, $\gd_2(\gep/2,2D+1,\eta/2)$,
$\gd_2(\gep/2,D,\eta)$, $\gd_1(\gep/2,2D+1,\eta/2)$. Assume that 
$\ga,-\gb\in(\eta,\pi-\eta)$ are such that there exists a curve
$f:[a,b]\to\BbC$ in $C_r(u,v)$ with $\diam(f)\leq D$, and let
$\bar\ga,\bar\gb$ be such that $|\bar \alpha -\alpha |,|\bar \gb -\gb |<\delta$ and
$C_r(\bar u, \bar v)$ is nonempty.
Although we cannot claim that both $C_r(\bar u, v)$ and
$C_r(u,\bar v)$ are nonempty (a counterexample can be easily found when
$f$ is a u-turn), it is easy to see that at least one of them is nonempty.
We will address the case when  $C_r(\bar u, v)$ is nonempty, as the other
case is similar. Since $|\bar \alpha -\alpha |<\gd_3\leq \gd_1(\gep/2,D,\eta)$,
we obtain from Lemma \ccurveone\ that there exists a curve $c_1\in C_r(\bar u,  v)$ such that
$||c_1||^2 \leq ||f||^2 + \gep/2$ and $\diam(c_1)\leq 2D+1$. Note that 
$\bar\ga,-\gb\in(\eta/2,\pi/\eta/2)$ (since $\gd_3\leq \eta/2$) and that $c_1$ is a
curve in $C_r(\bar u,  v)$ with $\diam(c_1)\leq 2D+1$. Since $C_r(\bar u, \bar v)$ is nonempty
and $|\bar \gb -\gb |<\gd_3\leq \gd_2(\gep/2,2D+1,\eta/2)$, it follows from Remark \ccurvetwo\ that
there exists a curve $c\in C_r(\bar u, \bar v)$ such that 
$\norm{c}^2< \norm{c_1}^2 + \gep/2<\norm{f}^2+\gep$.
\qed\enddemo

For given $\alpha ,\beta \in (-\pi, \pi)$ let $u=(0,e^{i\ga})$ and $v=(1,e^{i\gb})$ be the unit
tangent vectors with base points $0$ and $1$ and direction angles $\ga$ and $\gb$, respectively.
It will be useful to
introduce the following notations: $C_r(\alpha, \beta ):= C_r(u,v)$,
$S(\alpha, \beta ):= S(u,v)$ and $\Energy(\alpha,
\beta):=\Energy(u,v)$, which we will call the bending energy of
the pair $\alpha, \beta$. We say $\alpha, \beta $ are s-feasible
if $S(u,v)$ is nonempty.

\proclaim {Proposition \Tlabel \small  } The bending energy $\Energy(\ga,\gb)$
tends to $0$ as $\ga,\gb\to 0$.
\endproclaim

\demo {Proof} Assume that $\abs\ga,\abs\gb\leq\pi/3$, and let $f:[0,1]\to\BbC$ be the
curve given by $f(t)=t+ig(t)$, where
$g$ is the cubic polynomial $g(t)=t(\tan\ga-(\tan\ga+\tan\gb)t)(1-t)$. Then $f$ belongs
to $S(\ga,\gb)$ and it is easy to see that there exists a constant $C$ such that
$\dsize\norm{f}^2\leq \frac C4 \int_0^1|g''(t)|^2\dd t$. A simple calculation shows that
the latter quantity equals $C(\tan^2\ga+\tan\ga\,\tan\gb+\tan^2\gb)$ which tends to $0$
as $\ga,\gb\to 0$.
\qed\enddemo

In what follows we will rely heavily on the results of
sections 5 and 6, where it
is assumed that $u$ and $v$ are in `canonical'
arrangement: $\alpha \geq |\beta|$. However, if we perturb
$\alpha $ and $\beta$, the resultant pair, $\bar
\alpha$ and $\bar \beta$, may no longer be in `canonical'
arrangement. The following two propositions will help deal with
this situation.

\proclaim {Proposition \Tlabel \alfa }

(i) $\Energy(\alpha, \beta )=\Energy(-\beta,-\alpha)$

(ii) $\Energy(\alpha, \beta )=\Energy(-\alpha, -\beta )$

(iii) $\Energy(\alpha, \beta )=\Energy(\beta, \alpha )$.
\endproclaim

\demo {Proof} For any curve in $S(\alpha,
\beta )$ if we reflect the curve across the $x=1/2$ line and
reverse its orientation, we obtain a curve in $S(-\beta ,
-\alpha)$ with the same bending energy. This means that there is a
bijection between $S(\alpha, \beta )$ and $S(-\beta , -\alpha)$
which preserves the bending energy. This implies (i). Similarly,
reflection across the $x$-axis gives a bending energy-preserving
bijection between $S(\alpha, \beta )$ and $S(-\alpha, -\beta )$,
which yields (ii). Combining (i) and (ii) we obtain (iii).
\qed\enddemo

\proclaim {Proposition \Tlabel \alfabar } Let us assume that
$\alpha \geq |\beta|$ but $0<\bar \alpha <|\bar \beta|$.

Define $(\wt \alpha , \wt \beta):=\bigg \lbrace \eqalign {
(\bar \beta , \bar \alpha) \quad &\hbox {if}\quad 
\bar \beta>0 \cr (-\bar \beta, -\bar \alpha ) \quad &\hbox
{if}\quad \bar \beta<0 \cr }.$

If $|\alpha -\bar \alpha |, |\beta -\bar \beta |<\delta$, then
$|\alpha -\wt \alpha |, |\beta -\wt \beta |<\delta$ and $\wt\ga\geq|\wt\gb|$.
\endproclaim

\demo {Proof} The proof is elementary
and we will leave it to the reader. \qed\enddemo

Let us indicate how we will use the previous two propositions in
the proof of the next one. Assume that  $\alpha \geq |\beta|$, and
$|\alpha -\bar \alpha |, |\beta -\bar \beta |<\delta$ with $
\alpha >\delta >0$. If $\bar \alpha <|\bar \beta|$, then we will
replace $\bar \alpha, \bar \beta$ with a new pair $\wt \alpha,
\wt \beta $ as in Proposition {\alfabar}. Then the new pair
will be in `canonical' arrangement ($\wt \alpha \geq |\wt
\beta|$), $|\alpha -\wt \alpha |, |\beta -\wt \beta |<\delta
$ and from Proposition {\alfa} we have $\Energy(\alpha,
\beta)=\Energy(\wt \alpha, \wt \beta)$.

\proclaim {Proposition \Tlabel \scurve} With the notations
introduced above let us assume that  $\alpha , \beta $ are
s-feasible with $\alpha\geq |\beta|$, and let
$\eta\in(0,\pi/4)$.  If $\ga\in(\eta,\pi-\eta)$, then for every $\gep >0$ there is a 
$\gd=\delta(\gep, \eta)>0$ (depending only on $\gep, \eta$) such that
if $|\bar \alpha -\alpha|,|\bar \beta -\beta|<\delta $ and $\bar
\alpha, \bar \beta$ are s-feasible, then $\Energy(\bar \alpha,
\bar \beta)< \Energy(\alpha , \beta)+ \gep .$
\endproclaim

Before we start the proof of Proposition {\scurve} let us recall
some quantities defined in Definition {\Gdefine}.
For $\gamma \in \gG=[\alpha -\pi, \beta]\cap
(-\infty ,0)$ we have
$$\align
y_1(\alpha, \gamma )&={1\over 2}\int _0^{\alpha -\gamma }\sqrt
{\sin \tau } \dd\tau,\qquad
y_2(\beta, \gamma )={1\over 2}\int _0^{\beta -\gamma }\sqrt
{\sin \tau } \dd\tau,\\
G(\alpha, \beta, \gamma )&= {(y_1+y_2)^2\over -\sin \gamma },\qquad
G_{min}(\alpha, \beta)=\min \{G(\alpha, \beta,\gamma ):
\gamma \in \gG\}.
\endalign
$$
Recall from Section 3 that $d=\xi (\pi)={1\over 2}\int _0^{\pi
}\sqrt {\sin \tau } \dd\tau$ and define the quantity $\gamma _0$ by
$\gamma _0=-\sin ^{-1}\big((\sin \eta)(1-\cos \eta)^2/(16d^2)\big)$.
If $\ga\in(\eta,\pi-\eta)$, from the formulas above one can
verify immediately that
$$G(\alpha, \beta ,\alpha -\pi)\leq {4d^2\over \sin \eta}<
G(\alpha,\beta,\gamma)\quad \hbox {if} \quad \gamma _0<\gamma
<\beta. \tag \Elabel \gnul$$ This implies that if
$G(\alpha,\beta,\gamma)=G_{min}(\alpha, \beta)$, then $\gamma \leq
\beta^*=\min \{\beta,\gamma _0\}$.

It will be convenient to extend the domain of
$G(\alpha,\beta,\gamma)$ to include any $\gamma \in [-\pi, 0]$
without changing the minimum $G_{min}(\alpha, \beta)$ or values
$\gamma$ where the minimum is assumed.
We define the set $K_{\eta}$ by
$K_{\eta}=\{ (\alpha,\beta): \eta\leq \alpha \leq \pi-\eta,\quad
|\beta|\leq \alpha,\quad \alpha -\pi\leq \beta\}.$ For $(\alpha,
\beta)\in K_{\eta},\quad \beta ^*:=\min \{\beta,\gamma _0\}$ and
$\gamma\in [-\pi,0)$ we set
$$\widehat G(\alpha,\beta,\gamma)=\Bigg \lbrace \eqalign { &G(\alpha, \beta, \beta ^*)+\gamma -\beta
^*\cr
 &G(\alpha, \beta, \gamma )\cr  &G(\alpha, \beta, \alpha
-\pi)+\alpha-\pi-\gamma \cr  }\quad \eqalign { &\hbox {if}\quad
\beta ^*<\gamma \leq 0 \cr &\hbox {if}\quad \alpha -\pi\leq
\gamma\leq \beta ^* \cr &\hbox {if}\quad -\pi\leq \gamma < \alpha
-\pi \cr}  .$$ From the remark following inequality ({\gnul}) and
from the construction of $\widehat G$ it is clear that
$G_{min}(\alpha, \beta)=\widehat G_{min}(\alpha, \beta)=\min
\{\widehat G(\alpha, \beta,\gamma ): \gamma \in [-\pi , 0]\}$.
Moreover $G$ and $\widehat G$ assume their minimum at the same
points, that is  $G_{min}(\alpha, \beta)=G(\alpha,\beta,\gamma)$
if and only if $\widehat G_{min}(\alpha, \beta)=\widehat
G(\alpha,\beta,\gamma)$.

 The quantity $\sigma (\gamma)$ will be interesting for us
only in the case when $\gamma =\beta$. Therefore we have
$$ \sigma (\alpha , \beta)=\cos \beta +{\sin \beta \over
y_1(\alpha, \beta)}\sqrt {\sin (\alpha -\beta)},\qquad \alpha
-\pi\leq  \beta <0.$$

It is easy to see that $\sigma (\alpha , \beta )\to 1$ as $\beta
\to 0$. Therefore we can extend the domain of $\sigma (\alpha,
\beta )$ to the region $ 0\leq \beta $ by setting
$$\sigma (\alpha , \beta)=1, \qquad \hbox {if}\qquad \beta \geq 0.$$

\noindent We can summarize the results of sections 5 and 6 as follows:\newline

If $\sigma (\alpha, \beta )\leq 0$, then we are in Case C of Summary \SummarySuv\ 
and there is a segment of rectangular
elastica in $C_r(\ga,\gb)$ which has minimal bending energy in $S(\alpha, \beta )$.

If $\sigma (\alpha, \beta )>0$, then we are in Case A or B and there is an s-curve in
$S(\alpha, \beta )$ with minimal bending energy 
$\Energy(\alpha , \beta )=G_{min}(\alpha, \beta).$

\demo {Proof of Proposition  {\scurve }}

\noindent{\bf Case 1:} $\sigma (\alpha , \beta)\leq 0$. \newline
Then $\gb<0$ and there is a segment of rectangular
elastica $f\in C_r(\ga,\gb)$ which has minimal bending energy in $S(\alpha, \beta )$.
Before invoking Proposition \ccurve, we remark that the diameter of $f$ cannot exceed $10$
since the ratio of length over breadth for any segment of rectangular elastica is
bounded by $10$.
Since $\sigma (\alpha , \beta)$ is continuous and $\ga\in(\eta,\pi-\eta)$,
one can see that there is an $\eta_1=\eta _1(\eta) \in(0,\eta)$ such that
$|\beta| >\eta _1(\eta)$. By Proposition {\ccurve},
there is a $\gd=\delta (\eta, \gep)>0$ (namely, $\gd_3(\gep/2,10,\eta_1)$ in the language of
Prop. \ccurve) such that if 
$|\bar \alpha -\alpha|,|\bar \beta -\beta|<\delta $, then 
$$\Energy(\bar \alpha, \bar \beta)<\Energy(\alpha , \beta )+\gep/2,$$
which completes the proof for Case 1.

\noindent {\bf Case 2:} $\sigma (\alpha , \beta)> 0$. \newline
Then $\Energy(\alpha , \beta )=G_{min}(\alpha, \beta)=G(\alpha,
\beta,\gamma )$, for some (not necessarily unique) $\gamma =\gamma (\alpha, \beta)$. 
From the remark following inequality ({\gnul})
we have   $\gamma (\alpha, \beta)\in [\alpha -\pi, \beta ^*]$,
where $\beta ^*=\min \{\beta,\gamma _0\}$.
From the definition of $G$ and $\widehat G$ one can  see that 
$\widehat G(\alpha, \beta,\gamma )$ is continuous, hence uniformly
continuous on the region $K_{\eta}\times [-\pi, 0]$. Therefore,
there is a $\gd_0=\delta_0(\eta,\gamma _0)>0$ such that for $(\alpha,
\beta),( \alpha ',  \beta ')\in K_{\eta}$ we have
$$|\widehat G(\alpha, \beta,\gamma )-\widehat G( \alpha ',  \beta ',
\gamma' )|<{\gep \over 2},\qquad \hbox {whenever}\qquad |
\alpha-\alpha '|, | \beta-\beta '|,| \gamma-\gamma
'|<\delta_0.\tag \Elabel \gam
$$
Let us assume that $|\bar \alpha-\alpha |, |\bar
\beta-\beta|<\delta$, where $\delta \leq \min \{\delta _0, {\eta \over 4}\}$ is determined later. 
We can further assume, without loss of
generality, that $\bar \alpha \geq |\bar \beta|$, since otherwise, we can
replace $\bar \alpha, \bar \beta $ with $\wt \alpha, \wt \beta$, keeping in mind
that $|\wt \alpha-\alpha |, |\wt\beta-\beta|<\delta$, by Proposition {\alfabar},
and $\Energy(\wt \alpha, \wt \beta)=\Energy(\alpha, \beta)$,
by Proposition {\alfa}.

{\bf Case 2a:} $\sigma (\bar \alpha, \bar \beta )\geq 0$. \newline
Let $\bar \gamma = \gamma (\bar \alpha, \bar \beta)$ be an angle
where $G(\bar \alpha, \bar \beta,\gamma)$ assumes its minimum.
Since $G$ and $\widehat G$ assume their minimum at the same points
(see the remarks following the definition of $\widehat G$) we have
$ \Energy(\bar \alpha, \bar \beta )=G_{min}(\bar \alpha, \bar \beta )=G(\bar \alpha,
\bar \beta,\bar \gamma )= \widehat G(\bar \alpha, \bar \beta, \bar
\gamma )$ and
$\Energy( \alpha, \beta )=G_{min}( \alpha,  \beta )=G( \alpha,
 \beta, \gamma )= \widehat G( \alpha,  \beta,
\gamma )$.  Taking (\gam) into  consideration, we then obtain
$$ \Energy(\bar \alpha, \bar \beta )
=\widehat G(\bar \alpha,
\bar \beta,\bar \gamma )\leq \widehat G(\bar \alpha, \bar \beta,
\gamma )< \widehat G( \alpha,  \beta, \gamma )+{\gep \over
2}=\Energy(\alpha, \beta)+{\gep \over 2},$$
which completes the proof for Case 2a (with $\gd=\min \{\delta _0, {\eta \over 4}\}$).

{\bf Case 2b:} $\sigma (\bar \alpha, \bar \beta )<0$.\newline
Since
$\delta \leq \eta/4$ we have $(\alpha, \beta),(\bar \alpha, \bar
\beta)\in K_{\eta /2}$, which is a convex set. Therefore the line
segment  $[(\alpha, \beta),(\bar \alpha, \bar \beta)]$ is also a
subset of $K_{\eta /2}$. Since $\sigma$ is a continuous function
of $\alpha, \beta$ there is a pair $(\alpha _1, \beta_1)\in
[(\alpha, \beta),(\bar \alpha, \bar \beta)]$ with $\sigma (
\alpha_1,  \beta_1 )= 0$.
 Applying the previous argument for $\alpha _1, \beta_1$ instead
 of $\bar \alpha , \bar \beta$ we obtain
 $\Energy(\alpha _1,\beta _1)<\Energy(\alpha , \beta ) +{\gep \over
 2}.$

Since $\sigma ( \alpha_1,  \beta_1 )= 0$ there is a segment of
rectangular elastica $f\in C_r(\alpha_1,\beta_1)$  with $||f||^2=\Energy(\alpha _1,\beta
_1)$. Noting that $\alpha _1
>{\eta\over 2}$, we have, as in Case 1,  $|\beta
_1|>\eta _1(\eta/2)$. We can now apply Proposition {\ccurve} to obtain
$\Energy(\bar \alpha, \bar \beta)<\Energy( \alpha_1, \beta
_1)+{\gep \over 2},$\newline
provided $\delta\leq \delta_3 (\gep /2, 10, \eta _1(\eta/2))$ in the language of Proposition
{\ccurve}.
Combining this with the previous inequality we obtain
$$\Energy(\bar \alpha, \bar \beta)<\Energy(\alpha , \beta ) +{\gep },$$
which completes the proof for Case 2b. \qed\enddemo

With propositions \ccurve--\scurve\ in hand, we can finally prove the following.

\proclaim {Theorem \Tlabel \continuous } The bending energy
$\Energy(u,v)$ depends continuously on the directions of the unit tangent vectors $u$ and $v$.
\endproclaim

\demo {Proof} Let $u$ and $v$ be s-feasible unit tangent vectors. 
As explained at the beginning of section 5, we
can assume, without loss of generality,
that $u=(0,e^{i\ga})$ and $v=(1,e^{i\gb})$, where $\ga,\gb\in(-\pi,\pi)$ satisfy 
$\ga\geq\abs{\gb}$. The case $\ga=0$ has been settled
in Proposition \small, so assume $\ga\in(0,\pi)$. Let $\eta\in(0,\pi/8)$ be such that
$\ga\in(2\eta,\pi-2\eta)$. 

Let $\gep>0$ and set $\gd=\min\set{\eta,\gd_4}$, where
$\gd_4=\gd(\gep,\eta)$ is as described in  Proposition
\scurve. Let $\bar\ga,\bar\gb$ be s-feasible with  $|\bar \alpha -\alpha|, |\bar \beta -\beta|<\delta $.
We will show, in two cases, that
$$\abs{\Energy(\bar \alpha,\bar \beta)-\Energy( \alpha, \beta)}<\gep.
\tag \Elabel \driedtea$$

\noindent{\bf Case 1:} $\bar\ga\geq|\bar\gb|$.\newline
As written, Proposition \scurve\ yields  $\Energy(\bar \alpha,\bar \beta)< \Energy(\alpha , \beta)+\gep$,
but Proposition \scurve\ can also be applied with $(\ga,\gb)$ and $(\bar\ga,\bar\gb)$ interchanged,
since $\bar\ga\in(\eta,\pi-\eta)$ and $\bar\ga\geq|\bar\gb|$. This yields 
$\Energy(\alpha, \beta)<\Energy(\bar\alpha,\bar \beta)+\gep$, and we obtain (\driedtea).

\noindent{\bf Case 2:} $\bar\ga<|\bar\gb|$.\newline
Let $\wt \alpha, \wt\beta$ be as defined in Proposition {\alfabar}, whereby
$|\wt \alpha-\alpha |, |\wt \beta-\beta|<\delta$ and $\wt\ga\geq|\wt\gb|$. It follows from Case 1 that
$\abs{\Energy(\wt \alpha,\wt \beta)-\Energy( \alpha, \beta)}<\gep$, but
since  $\Energy(\wt \alpha, \wt \beta)=\Energy(\bar\alpha, \bar\beta)$ (by Proposition \alfa),
we have (\driedtea).
\qed\enddemo

\demo {Proof of  Theorem \theoremmain } Most of the proof has been explained in the discussion
following the theorem's statement at the beginning of this section.  All that remains is to show that
our candidate $c=s_1\cup s_2 \cup \cdots \cup s_{m-1}$ has bending energy $M$ 
(the infimum of bending energies in $\Cal A(P_1,P_2,\ldots,P_m)$).

For $n=1,2,3,\ldots$, let us write $c^n \in \Cal
A(P_1,\ldots ,P_m)$ as $c^n=f_1^n\splice f_2^n \splice \cdots \splice f_m^n$, where $f_j^n$ belongs to $S(v_j^n,v_{j+1}^n)$
for $j=1,2,\ldots,m$. Note that
$$\norm{c^n}^2 = \sum_{j=1}^{m-1} \norm{f_j^n}^2 \geq \sum_{j=1}^{m-1} \Energy(v_j^n,v_{j+1}^n)\quad\text{ while }\quad
\norm{c}^2 = \sum_{j=1}^{m-1} \norm{s_j}^2 = \sum_{j=1}^{m-1} \Energy(v_j,v_{j+1})$$ 
(since $s_j$ has minimal bending energy in $S(v_j,v_{j+1})$).
For $j=1,2,\ldots,m-1$, it follows from Theorem \continuous\ that 
$\Energy(v_j^n,v_{j+1}^n)\to \Energy(v_j,v_{j+1})$ as $n\to\infty$, and therefore we obtain
$\norm{c}^2 \leq \lim_n \norm{c^n}^2 = M$. Since $c$ belongs to $\Cal A(P_1,P_2,\ldots,P_m)$, 
we conclude that $\norm{c}^2=M$.
\qed\enddemo

\remark{Remark \Tlabel\remarrk} We conjecture, and hope to show in a subsequent
paper, that if each of the curves $s_i$ is of form one (see
Definition {\firstsecond}), then the resulting optimal curve $s$
is twice continuously differentiable.
\endremark

\remark{Remark \Tlabel\rremarrk} Let us denote by $\Cal
A_{periodic}(P_1,P_2,\ldots ,P_m,P_1)$ the set of periodic (closed) curves passing
through the points $P_1,P_2,\ldots ,P_m,P_1$ such that they are s-curves
between any two consecutive points. Notice that the proof of
Theorem {\theoremmain} works equally well for periodic
admissible curves, provided that $\Cal A_{periodic}(P_1,P_2,\ldots ,P_m,P_1)$ is nonempty. 
This is exactly what is shown
in the proof of Proposition \nemures. Therefore we have the
following extension of Theorem {\theoremmain }:
\endremark

\proclaim {Theorem \Tlabel\maintt  } Given any sequence of points
$P_1,P_2,\ldots ,P_m\in \BbC$ with $P_j\ne P_{j+1}$ and $P_m\ne P_1$, the family
$\Cal A_{periodic}(P_1,P_2,\ldots ,P_m,P_1)$ contains a curve with minimal
bending energy.
\endproclaim

\subhead Acknowledgements \endsubhead
The authors are very grateful to Hakim Johnson (Kuwait English School) for writing the computer program 
{\it Curve Ensemble}, based on elastic splines, which was used to make the figures. We are also
grateful to Aurelian Bejancu for discussions on variational calculus which led to a clean proof of Theorem 3.2,
and to the referees and editor for many helpful comments and suggestions.
%

%
\enddocument